\documentclass{article}%
\usepackage{amsmath}
\usepackage{amsfonts}
\usepackage{amssymb}
\usepackage{graphicx}
\usepackage[pdftex]{geometry}
\usepackage[pdftex]{color}
\usepackage{tikz,pgfplots}
\usepackage{subcaption}%
\providecommand{\U}[1]{\protect\rule{.1in}{.1in}}
\newtheorem{theorem}{Theorem}

\newtheorem{conjecture}[theorem]{Conjecture}
\newtheorem{corollary}[theorem]{Corollary}

\newtheorem{lemma}[theorem]{Lemma}

\newtheorem{proposition}[theorem]{Proposition}
\newtheorem{remark}[theorem]{Remark}

\newenvironment{proof}[1][Proof]{\noindent\textbf{#1.} }{\ \rule{0.5em}{0.5em}}
\pgfplotsset{compat=1.12}
\begin{document}

\title{Stability estimate for the {H}elmholtz equation with rapidly jumping
coefficients \thanks{The authors gratefully acknowledge support by the Swiss
National Science Foundation under grant no. 172803. The authors are also
grateful to the Applied Mathematics Department at ENSTA ParisTech for the kind
hospitality in the fall semester 2016 and the Hausdorff Research Institute for
Mathematics in Bonn for Visiting Fellowships in their 2017 Trimester Programme
on Multiscale Methods, during which part of this work was carried out.}}
\author{Stefan Sauter\thanks{(stas@math.uzh.ch), Institut f\"{u}r Mathematik,
Universit\"{a}t Z\"{u}rich, Winterthurerstr 190, CH-8057 Z\"{u}rich,
Switzerland}
\and Celine Torres \thanks{(celine.torres@math.uzh.ch), Institut f\"{u}r
Mathematik, Universit\"{a}t Z\"{u}rich, Winterthurerstr 190, CH-8057
Z\"{u}rich, Switzerland}}
\maketitle

\begin{abstract}
The goal of this paper is to investigate the stability of the Helmholtz
equation in the high-frequency regime with non-smooth and rapidly oscillating
coefficients on bounded domains. Existence and uniqueness of the problem can
be proved using the unique continuation principle in Fredholm's alternative.
However, this approach does not give directly a coefficient-explicit energy
estimate. We present a new theoretical approach for the one-dimensional
problem and find that for a new class of coefficients, including coefficients
with an arbitrary number of discontinuities, the stability constant (i.e., the
norm of the solution operator) is bounded by a term independent of the number
of jumps. We emphasize that no periodicity of the coefficients is required. By
selecting the wave speed function in a certain \textquotedblleft
resonant\textquotedblright\ way, we construct a class of oscillatory
configurations, such that the stability constant grows exponentially in the
frequency. This shows that our estimates are sharp.

\end{abstract}



\section{Introduction\label{SecIntro}}

\label{Sec:intro}

High-frequency scattering problems have many important applications which
include, e.g., radar and sonar detection as well as medical and seismic
imaging. In physics, such problems are studied intensively in the context of
wave scattering in disordered media and localization of waves with the goal to
\textit{design} waves with prescribed intensity, interference, localized foci,
parity-time symmetry, etc. Important applications are in nano photonics and
lasers -- see, e.g., \cite{Anderson_local}, \cite{LAGENDIJK1996143},
\cite{sebbah2001waves}, \cite{lahini2008anderson}, \cite{makris2017wave},
\cite{transmission_matrix_dis_sys}, \cite{lowe1995matrix} for references to
the theoretical and experimental physics literature.

Their efficient and reliable numerical modelling is a challenge and the
development of fast numerical methods is far from being mature. Such problems
are often modelled in the frequency domain, where a time-periodic ansatz is
employed for the wave equation which, in turn, results in the Helmholtz
equation in the high-frequency regime, i.e., with large wave number. Moreover,
in applications such as seismic or medical imaging, the media typically are
heterogeneous and, consequently, the coefficients in the Helmholtz equation
become variable. The numerical analysis for these types of problems is much
less developed as the high-frequency homogeneous case.

In this paper, we discuss the Helmholtz equation of the form%
%TCIMACRO{\TeXButton{HelmStrongForm}{\begin{subequations}
%\label{HelmStrongForm}
%\end{subequations}}}%
%BeginExpansion
\begin{subequations}
\label{HelmStrongForm}
\end{subequations}%
%EndExpansion%
\begin{equation}
-\mathrm{div}(a~\mathrm{grad}~u)-\left(  \frac{\omega}{c}\right)
^{2}u=f\text{ in }\Omega,\tag{%
%TCIMACRO{\TeXButton{HelmStrongForm}{\ref{HelmStrongForm}}}%
%BeginExpansion
\ref{HelmStrongForm}%
%EndExpansion
a}\label{HelmStrongForma}%
\end{equation}
on a bounded Lipschitz domain $\Omega$, frequency $\omega>0$, positive wave
speed $c$ and diffusion coefficient $a$, where both, $c$ and $a$, are
variable. On the Helmholtz problem, we impose impedance boundary conditions%
\begin{equation}
a\frac{\partial u}{\partial\mathbf{n}}-\operatorname*{i}\sqrt{a}\frac{\omega
}{c}u=g\text{ on }\Gamma=\partial\Omega.\tag{%
%TCIMACRO{\TeXButton{HelmStrongForm}{\ref{HelmStrongForm}}}%
%BeginExpansion
\ref{HelmStrongForm}%
%EndExpansion
b}\label{HelmStrongFormb}%
\end{equation}
Let \(L^2(\Omega)\) be the usual Lebesgue space with scalar product and norm \[(u,v):= \int_\Omega u \bar{v},\qquad \|u\|:= (u,u)^2.\] We define the \textquotedblleft energy space\textquotedblright\ by
$\mathcal{H}:=H^{1}\left(  \Omega\right)  $, equipped with
\begin{align}\label{def:energynorm}
\left\Vert w\right\Vert _{\mathcal{H}}:=\sqrt{\left\Vert \sqrt{a}\nabla
w\right\Vert ^{2}+\left\Vert \frac{\omega}{c}w\right\Vert ^{2}}.
\end{align}
The
variational formulation of (\ref{HelmStrongForm}) is to find $u\in\mathcal{H}$
such that%
\begin{equation}
B\left(  u,v\right)  :=\left(  a\nabla u,\nabla v\right)  -\left(
\frac{\omega}{c}u,\frac{\omega}{c}v\right)  -\operatorname*{i}\left(  \sqrt
{a}\frac{\omega}{c}u,v\right)  _{L^{2}\left(  \Gamma\right)  }=\left(
f,v\right)  +\left(  g,v\right)  _{L^{2}\left(  \Gamma\right)  }=:F\left(
v\right)  \quad\forall v\in\mathcal{H}.\label{HelmVarForm}%
\end{equation}

To solve the problem numerically one may use abstract, conforming Galerkin
methods, i.e., choose a finite-dimensional subspace $S\subset\mathcal{H}$ and
seek for $u_{S}\in S$ such that%
\begin{equation}
B\left(  u_{S},v\right)  =F\left(  v\right)  \qquad\forall v\in S\text{.}
\label{GalDiscSNF}%
\end{equation}

The stability constant $C_{\operatorname{stab}}$ of the problem, satisfying
\begin{align*}
\|u\|_{\mathcal{H}} \leq C_{\operatorname{stab}} \left(  \|f\|^{2}+\|g\|_{H^{1/2}(\Gamma)}^{2} \right)  ^{\frac{1}{2}},
\end{align*}
plays an important role for the design and the numerical analysis of abstract
Galerkin methods of the form (\ref{GalDiscSNF}) and much research has been
devoted to estimate this constant. $C_{\operatorname*{stab}}$ in general
depends on the wave number and the coefficient functions $a$ and $c$.

In two dimensions and for $a$, $c$ positive constants, $\Omega$ convex or
$C^{1}$-star shaped the stability constant $C_{\operatorname{stab}}>0$ is
independent of $f$, $g$, and $\omega$ \cite{MelenkDiss}. The result was
extended to higher dimensions in \cite{CummingsFeng2006}. For the case of
general Lipschitz domains and constant coefficients, it was proved in
\cite[Theorem 2.4]{Esterhazy2012} that $C_{\operatorname*{stab}}\leq
C\omega^{5/2}$ while the result was improved in \cite{Spence2014} to
$C_{\operatorname*{stab}}\leq C\omega$. For bounded domains with $C^{\infty}$
boundaries it was proved in \cite[Theorem 1.8]{baskin2016sharp} that
$C_{\operatorname*{stab}}=O\left(  1\right)  $. The convergence analysis for
$hp$-Galerkin finite elements as in \cite{MelenkSauter2010,MelenkSauter2011}
shows that in the error estimate $C_{\operatorname*{stab}}$ is multiplied by a
term, which is exponentially small in the polynomial degree $p$, so that a
choice $p\sim\log(\omega)$ preserves the optimal convergence order without pollution.

\paragraph{Stability Estimates for the heterogeneous Helmholtz Problem.}

The theory of the heterogeneous Helmholtz equation is much less developed than
the one for constant coefficients since the wave strongly depends on
interfaces, variations of the coefficients, and can exhibit localization,
interferences, complicated variations of intensity, in particular, if
scattering through disordered media is considered, see, e.g.,
\cite{sebbah2001waves}, \cite{MoiolaSpence_resonance_2017}.

First results regarding rigorous stability estimates in one dimension go back
to \cite{AKS1988}, where $c$ is supposed to be sufficiently smooth (at least
$C^{1}$). For $d\geq2$ and slowly varying coefficients, similar results were
proved in \cite{PerthameVega1999} and \cite{GrahamSauter2017}. The main
theoretical approaches for estimating the stability in the case of varying
coefficients are as follows.

a) \textquotedblleft Rellich\textquotedblright\ or \textquotedblleft
Morawetz\textquotedblright\ multipliers. For general dimension a test function of the form
$v=x.\nabla u$ and $v=u$, or modifications thereof is employed in the
variational formulation of the Helmholtz equation which allows to estimate the
$L^{2}$-part of the solution in terms of its $H^{1}$-semi norm and the data.
In turn, this estimate is used to estimate the stability constant. This technique was further developped for the heterogeneous Helmholtz equation in \cite{PerthameVega1999}.

In \cite{ChaumontFreletDiss, GrahamSauter2017}, stability is achieved via the
\textquotedblleft Rellich trick\textquotedblright, by replacing the smooth
field $x$ with a coefficient-dependent piecewise smooth field. It turns out that,
only oscillatory wave speeds which are small perturbations of a constant wave
speed, can be handled by this technique .

b) Full space methods. If heterogeneous Helmholtz problems are considered in
full space, methods from semi-classical/asymptotic/microlocal analysis can be
applied and estimates for the stability constant are derived, e.g., in
\cite{Burq_hetero_analysis} for smoothly varying coefficients and in
\cite{Bellassoued_carleman} for a full space problem with one inclusion and a
discontinuity across the interface. It is shown that for these cases, the
stability constant can grow \textit{at most} exponentially in $\omega$. In
\cite{MR0457962}, and \cite{popov1999resonances},
\cite{MoiolaSpence_resonance_2017}, examples are presented for smooth
coefficients/discontinuous coefficients with one interface where the stability
constant grows super-algebraically.

c) Homogenization. For \textit{periodic}, heterogeneous media, methods of
homogenization can be applied (for diameters of inclusions tend to zero and
their number goes to infinity while the frequency $\omega$ is assumed to be
fixed) to derive effective equations (see, e.g.,
\cite{bouchitte2004homogenization}, \cite{ohlberger_verfuerth}) which then can
be analysed by methods for Helmholtz equations with constant coefficients.

d) Matrix techniques. In layered materials, one can approximate the wave
propagation problem by restricting to a linear combination of
\textit{finitely} many waves types (longitudinal and shear waves) in each
homogeneous material part of the domain as an ansatz and combine this across
the interfaces by transmission and boundary conditions (see, e.g.,
\cite{lowe1995matrix} and citations therein). This leads to a linear algebraic
system for the coefficients in this ansatz. Our approach is following the same
basic idea. However, since we are in 1D the fundamental system of the
Helmholtz equation consists of only two types of waves (one is outgoing, the
other one ingoing) and we end up with an $n\times n$ tridiagonal block system
consisting of $2\times2$ matrices per block -- here $n$ denotes the number of
discontinuities. Our main achievement in this paper is the derivation of a
stability bound of the corresponding Green's function which is explicit with
respect to the wave number, the coefficients, and the number of jumps. In the
case of a stochastic 1D medium the Green's function has been analysed in
\cite{frisch1973stochastic} and asymptotic properties have been derived in a
stochastic setting. Also in this case, realizations of the random coefficients
are considered where the wave exhibits localization and near-resonances.

\paragraph{Explicit Estimates in 1D for the Heterogeneous Helmholtz Problem.}

In contrast to the Green's function (of the ordinary differential equation)
for the Helmholtz problem in one dimension with constant coefficients, the
Green's function in the case of variable coefficients is not known explicitly
and more subtle mathematical tools have to be employed for 1D stability
estimates. In \cite{AKS1988}, a one-dimensional model problem for the
Helmholtz problem was considered and well-posedness of (\ref{HelmStrongForm})
was proved for positive and bounded wave speed $c$. A test function of
Rellich/Morawetz-type (cf. \cite{Rellich:40}, \cite{MorwaetzLudwig}) of the
form $v=au^{\prime}+bu$ for some functions $a$ and $b$ (chosen as solutions of
a certain ODE) has been employed to prove the wave-number independent
stability bounds for $c\in C^{1}\left(  \Omega\right)  $ while the constant
depends on $\left\Vert c^{\prime}\right\Vert _{\infty}$. In \cite{MaIhBa:96},
\cite{Ihlenburgbook}, the test function $v\left(  x\right)  =xu^{\prime
}\left(  x\right)  $ of Rellich/Morawetz-type has been used to prove stability
bounds for a certain fluid-solid interaction problem for elastic waves which
can be considered as a problem with a wave speed which has discontinuities at
two points. For the one-dimensional problem, some explicit estimates of the
stability constant were published only recently. In the thesis \cite{ChaumontFreletDiss},
 a stability estimate for the case of,
possibly, very large numbers of layers and general, layer-wise constant wave
speed $c$ and $a=1$ is proved. The estimates are explicit in the number of
jumps and the values of the coefficients. In the worst case, the estimate
grows exponentially in the number of interfaces. A generalized result is
proved in \cite{GrahamSauter2017}, where $a$, $c$ are both piecewise monotonic
functions with bounded variation. In this paper, we present estimates of the
stability constant in 1D for certain classes of piecewise constant wave speeds
$c\in L^{\infty}(\Omega,\mathbb{R}_{>0})$ which are explicit in the wave
number $\omega$ and $c$. We do not require any periodicity on $c$. This result
is finer than the result in \cite{GrahamSauter2017}, where the stability
constant grows exponentially in the variation of $c$ and $a$. In addition, we
will construct a class of configurations which show that a) our estimates are
sharp and b) that these configurations are very rare.

\paragraph{Outline of the Paper.}
Since the proofs in this paper are quite technical we explain our approach and
the main results here as an outline and summary. In Section
\ref{Sec:GeneralmultiD}, we discuss briefly the well-posedness of the problem
in the multidimensional case and state a conjecture on the behaviour of the
stability constant $C_{\operatorname{stab}}$ for a wave speed $c\in L^{\infty
}(\Omega,\mathbb{R}_{>0})$. We establish results towards this conjecture in 1D
in the following sections. In Section \ref{Sec:General1D}, we consider the
Helmholtz equation (with piecewise constant wave speed $c$, inhomogeneous
impedance boundary conditions, and without volume forces) and employ an ansatz
as a linear combination of ingoing and outgoing waves on each subinterval
combined by transmission and boundary conditions. This leads to a linear
$n\times n$ tridiagonal block system for the coefficients in this ansatz. The
blocks are symmetric (but not Hermitian)\ $2\times2$ matrices and $n$ denotes
the number of jumps of the wave speed. We derive a simple representation (see
Corollary \ref{CorrepR} for details) of the corresponding Green's function:%
\[
\left\vert \left(  \mathbf{M}_{\mathrm{Green}}^{\left(  2n\right)  }\right)
_{2n-2m+1,2n}\right\vert =\left\vert {\prod_{\ell=n-m+1}^{n}}\frac
{\sqrt{1-q_{\ell}^{2}}}{\left(  1+q_{\ell}Q_{\ell-1}\right)  }\right\vert
,\qquad1\leq m\leq n,
\]
where $q_{\ell}=\frac{c_{\ell+1}-c_{\ell}}{c_{\ell+1}+c_{\ell}}$ denote the
relative jumps on neighbouring intervals and $\left(  Q_{\ell}\right)  $ is a
recursively defined sequence with values in the interior of the complex unit
disc. This representation is based on Cramer's rule and a recursive formula
for the arising determinants. The proofs of these representations are very
technical but elementary and we have shifted them to the last section
(\S \ref{ProofsRep}). This representation allows us to reduce the estimate of
the solution operator to an analysis of the sequence $Q_{\ell}$, depending on
the piecewise constant values of the wave speed $c$ and a phase factor
$\sigma_{\ell}$ which encodes the\ interplay of the wave number $\omega$, the
values of the wave speed $c$, and the width of the subinterval $\tau_{\ell}$.
An estimate of the entries of Green's function is given in Lemma \ref{LemGnm}%
\[
\left\vert \left(  \mathbf{M}_{\mathrm{Green}}^{\left(  2n\right)  }\right)
_{2n-2m+1,2n}\right\vert \leq\frac{1}{\sqrt{1-\left\vert Q_{n-m}\right\vert
^{2}}},\qquad1\leq m\leq n
\]
and we recall $\left\vert Q_{\ell}\right\vert <1$. From this, we derive our
main stability estimate (Theorem \ref{Theo:StabEst})%
\begin{equation}
\left(  \int_{\Omega}|u^{\prime}|^{2}+\left(  \frac{\omega}{c}\right)
^{2}|u|^{2}\right)  ^{\frac{1}{2}}\leq C_{\operatorname{stab}}\max\left\{
\left\vert g_{1}\right\vert ,\left\vert g_{2}\right\vert \right\}
\quad\text{with\quad}C_{\operatorname*{stab}}:=4\frac{c_{\max}}{c_{\min}}%
\max_{1\leq\ell\leq n}\frac{1}{\sqrt{1-\left\vert Q_{\ell}\right\vert ^{2}}%
}.\label{Cstabintro}%
\end{equation}
This shows the importance of the sequence $\left(  Q_{\ell}\right)  _{\ell}$.
However, this sequence depends on a very high-dimensional parameter space
involving the location of the $n$ jump points, the wave number $\omega$,
and the values of the wave speed. Our focus is on large number of jumps $n$
and hence, the analysis of $Q_{\ell}$ becomes quite involved.

Before proving general estimates on $Q_{\ell}$ we present two classes of
examples in Section \ref{Sec:ex}. Both examples have in common that the number
of jumps tends to infinity and known \textit{a priori bounds} for the
stability constant (see, e.g., \cite{ChaumontFreletDiss},
\cite{GrahamSauter2017}) tend to infinity for these cases. However, the
near-resonance cases seem to be very rare and their construction require a
subtle tuning of the parameters (\S \ref{example_exp}) so that $\left\vert
Q_{\ell}\right\vert $ approaches, exponentially fast, the value $1$ and, in
turn, $C_{\operatorname*{stab}}$ grows exponentially with respect to $\omega$
(cf. (\ref{Cstabintro}), Remark \ref{Remlowerbound}). The construction of
\textit{well-behaved} examples is much simpler (\S \ref{example_osc});
although the bounds in \cite{ChaumontFreletDiss}, \cite{GrahamSauter2017} tend
to infinity as $n\rightarrow\infty$, our new bound (and of course the
solution) stays bounded independent of $n$. The implications are two-fold: the
results in \cite{ChaumontFreletDiss}, \cite{GrahamSauter2017} are sharp in
general, but very pessimistic in many parameter configurations.

Section \ref{Sec:inverseGreen} is devoted to the parameter-explicit estimate
of the sequence $\left(  Q_{\ell}\right)  _{\ell}$. In Lemma \ref{estA1Q1} we
employ tools from complex analysis to maximize the sequence with respect to
all parameters and derive%
\begin{equation}
\left\vert Q_{\ell-1}\right\vert \leq\frac{1-\kappa^{-\ell}}{1+\kappa^{-\ell}%
}\quad\text{with the \textit{condition number} of the wave speed\textit{ }%
}\kappa:=\left\Vert c\right\Vert _{\max}\left\Vert c^{-1}\right\Vert _{\max
}.\label{condnumber_wavespeed}%
\end{equation}
Hence, if the number of jumps are bounded from above by $O\left(
\omega\right)  $, the stability constant in (\ref{Cstabintro}) grows at most
exponentially in $\omega$ (see Section \ref{hbounded}). The analysis of the
remaining case, i.e., when the number of jumps is much larger than $\omega$,
is more involved and we restrict to a slightly simplified parameter setting in
Section \ref{hnotbounded}. We still allow the jump points to be distributed in
an arbitrary way while we assume that the wave speed is jumping between only
two positive values. Then it is again possible (Prop. \ref{Prop:hsmall}) to
maximize $\left\vert Q_{\ell}\right\vert $ for this setting. By employing the
side condition that the lengths of the subinterval with constant wave speed
have to sum up to the total length of $\Omega$ we derive an estimate of
$\left\vert Q_{\ell}\right\vert $ of the form%
\[
\left\vert Q_{\ell}\right\vert \leq1-C\alpha^{\omega/c_{\min}}%
\]
for some $\alpha\in\left(  0,1\right)  $. In combination with
(\ref{Cstabintro}) the conjectured bound in the considered scenarios follows
(\S \ref{SecStabEst}).

\section{High-Frequency Helmholtz Equations with Variable Coefficients}

\label{Sec:GeneralmultiD}

\subsection{Helmholtz Equation for Varying Coefficients}

We consider the Helmholtz Equation on a bounded Lipschitz domain
$\Omega\subset\mathbb{R}^{d}$ with variable wave speed $c$ and diffusion
coefficient $a$
\begin{equation}
-\mathrm{div}(a~\mathrm{grad}~u)-\left(  \frac{\omega}{c}\right)
^{2}u=f\text{ in }\Omega.\label{MultiDHelmholtz}%
\end{equation}
The right-hand side $f$ is in $L^{2}(\Omega)$ and we denote by $\omega$ the
frequency parameter, bounded from below by $\omega_{0}>0$. We assume that the
wave speed as well as the diffusion coefficient are bounded in the following
way
\begin{align*}
0<c_{\min}\leq &  c\leq c_{\max}<\infty,\\
0<a_{\min}\leq &  a\leq a_{\max}<\infty.
\end{align*}
We denote the boundary of $\Omega$ by $\Gamma:=\partial\Omega$ and let
$\Gamma_{\operatorname{R}},~\Gamma_{\operatorname{D}}$ be relatively open
pairwise disjoint subsets of $\Gamma$, with $\Gamma=\overline{\Gamma
_{\operatorname{R}}\cup\Gamma_{\operatorname{D}}}$. On (\ref{MultiDHelmholtz}%
), we impose the impedance boundary conditions
\begin{equation}
a\frac{\partial u}{\partial\mathbf{n}}-\operatorname*{i}\sqrt{a}\beta
u=g\text{ on }\Gamma_{\operatorname*{R}},\quad u=0\text{ on }\Gamma
_{\operatorname*{D}},\label{MultiDBoundaryConditions}%
\end{equation}
for a given impedance coefficient $\beta\in L^{\infty}(\Gamma
_{\mathrm{\operatorname*{R}}},[0,\beta_{\max}])$, $\beta_{\max}<\infty$. Let
$\mathcal{H}:=\{u\in H^{1}(\Omega)~|~u|_{\Gamma_{D}}=0\}$, with norm defined in \eqref{def:energynorm}.
The weak formulation of the problem is: For given
$\mathcal{F}\in\mathcal{H}^{\prime}$, find $u\in\mathcal{H}$ such that
\begin{equation}
B(u,v)=\mathcal{F}(v)\quad\forall v\in\mathcal{H},\label{MultiDProblem}
\end{equation}
where \(B(\cdot,\cdot)\) is defined as in \eqref{HelmVarForm}.

\subsection{Well-Posedness}

\begin{theorem}
Let $\beta\in L^{\infty}\left(  \Gamma_{\mathrm{\operatorname*{R}}}%
,\mathbb{R}_{\geq0}\right)  $ be such that $\operatorname*{supp}\beta$ has
positive $\left(  d-1\right)  $-dimensional boundary measure. Let $a,c\in
L^{\infty}(\Omega)$ with $0<a_{\min}\leq a^{\ast}\leq a_{\max}<\infty$. For
$d\geq3$ we assume, in addition, that $a\in C^{0,1}(\Omega)$. Then, for any
$\omega\geq\omega_{0}$, the heterogeneous Helmholtz problem
(\ref{MultiDProblem}) is well posed.
\end{theorem}

For $d=2$ the proof can be found in \cite[Theorem 1.1]{Alessandrini2012}, for
$d=1$ the proof is very similar. For $d=3$, one uses the fact that $a\in
C^{0,1}(\Omega)$ can be extended to $\mathbb{R}^{d}$. The proof of the theorem
is based on Fredholm's alternative (cf.
\cite{Alessandrini2012,GrahamSauter2017,JerisonKenig85}). This technique,
however, does only provide an implicit stability estimate. It is therefore not
straightforward, how the stability constant depends explicitly on $\omega$,
$c$ or other parameters.

\subsection{Maximal Growth of the Stability Constant}

As discussed, there has been several contributions for finding estimates of
the stability constant that are explicit in the parameters $\omega,a$ and $c$.
Recent results (cf. \cite{ChaumontFreletDiss, GrahamSauter2017}) show that the
stability constant can be bounded with respect to the number of jumps or the
total variation of the coefficients, if they are piecewise constant.
Considering coefficients that are highly oscillatory and have an increasing
number of jumps, the stability constant known from
\cite{ChaumontFreletDiss,GrahamSauter2017} are diverging to infinity. Based on
known theoretical results and also numerical experiments, we formulate the
following conjecture.

\begin{conjecture}
For any bounded Lipschitz domain $\Omega\subset\mathbb{R}^{d}$, $a=1$, $c\in
L^{\infty}(\Omega)$, with $0<c_{\min}\leq c\leq c_{\max}<\infty$, $\omega
\geq\omega_{0}>0$ it holds
\begin{equation}
\left(  \int_{\Omega}|\nabla u|^{2}+\left(  \frac{\omega}{c}\right)
^{2}|u|^{2}\right)  ^{\frac{1}{2}}\leq C_{\operatorname{stab}}\left(  \Vert
f\Vert_{L^{2}(\Omega)}^{2}+\Vert g\Vert_{L^{2}(\Gamma)}^{2}\right)  ^{\frac
{1}{2}}, \label{conjecture}%
\end{equation}
with
\[
C_{\operatorname{stab}}\leq C_{1}\exp\left(  C_{2}\omega\right)  ,
\]
$C_{1},C_{2}>0$ depending on $c_{\operatorname{min}},~c_{\operatorname{max}}$
and $\Omega$.
\end{conjecture}

For constant wave speed $c$, there is a fairly good knowledge about the
constant $C_{\operatorname*{stab}}$; under moderate assumptions on the domain
it can be shown that $C_{\operatorname*{stab}}$ only grows algebraically with
respect to $\omega$ and in many cases this constant is bounded independent of
$\omega$ (see the literature review in \S \ref{SecIntro}).

In higher dimension, there exist much less results. As explained in the
literature review (see \S \ref{SecIntro}), examples are known where
$C_{\operatorname*{stab}}$ grows at \textit{least} super-algebraically with respect to
$\omega$ while for smoothly varying coefficients and coefficients with a
discontinuity across one interface, an upper bound can be proved which grows
\textit{at most} exponentially in $\omega$. All these results are in accordance with
our conjecture.

In this paper, we present further results towards this conjecture for a
one-dimensional example. We will derive a recursive representation of the
Green's function which allows to understand well-behaved parameter
configurations as well as to detect near-resonance cases. In turn, this allows
us to find classes of wave speeds, where the stability constant grows
exponentially with respect to $\omega$.

\section{Helmholtz Equation in one Dimension}

\label{Sec:General1D} We subdivide the domain $\Omega=\left(  -1,1\right)  $
into $(n+1)$ intervals by introducing the mesh points%
\begin{equation}
-1=x_{0}<x_{1}<\ldots x_{n+1}=1\label{defxi}%
\end{equation}
and define the subintervals $\tau_{j}=\left(  x_{j-1},x_{j}\right)  $, $1\leq
j\leq n+1$ and widths%
\[
h_{j}:=x_{j}-x_{j-1}\text{.}%
\]
In this section, we consider piecewise constant wave speed which is given by%
\begin{equation}
c\left(  x\right)  :=c_{j}\text{ for }x\in\tau_{j},\label{defc}%
\end{equation}
where $c_{j}$ are positive constants. For a positive wave number $\omega
\in\mathbb{R}_{\geq\omega_{0}}$ we consider the homogeneous Helmholtz equation
in the strong form%
\begin{equation}%
\begin{array}
[c]{cl}%
-u^{\prime\prime}-\left(  \frac{\omega}{c}\right)  ^{2}u=0 & \text{in }%
\Omega=\left(  -1,1\right)  ,\\
-u^{\prime}-\operatorname*{i}\frac{\omega}{c_{1}}u=g_{1} & \text{at }x=-1,\\
u^{\prime}-\operatorname*{i}\frac{\omega}{c_{n}}u=g_{2} & \text{at }x=1.
\end{array}
\label{Helmpwconst}%
\end{equation}
The physical parameters of the problem are a) the jump points $\left(
x_{\ell}\right)  _{\ell=1}^{n}$ satisfying (\ref{defxi}), b) the (positive)
values $c_{\ell}$, $1\leq\ell\leq n+1$, of the wave speed on the subintervals,
c) the positive wavenumber $\omega>0$, and d) the right-hand side $g_{1}$,
$g_{2}$ in (\ref{Helmpwconst}). For the \textit{mathematical} analysis of the
problem, we introduce derived parameters (which will simplify the notation in
the proofs), namely, phase factors $\sigma_{j}\in\mathcal{C}:=\left\{
z\in\mathbb{C}:\left\vert z\right\vert =1\right\}  $ and relative jumps
$q_{j}$
\begin{subequations}
\label{defalphaij}%
\begin{align}
\sigma_{j}:=\exp\left(  -2\operatorname*{i}\frac{h_{j+1}\omega}{c_{j+1}%
}\right)  , &  \qquad0\leq j\leq n,\\
q_{j}:=\frac{c_{j+1}-c_{j}}{c_{j+1}+c_{j}}, &  \qquad1\leq j\leq n.
\end{align}

\end{subequations}
\begin{remark}
Note that $q_{j}\in\lbrack-q,q]$ holds for some $0<q<1$ and we have the
\textquotedblleft inverse\textquotedblright\ relation
\begin{equation}
c_{j+1}=c_{j}\frac{1+q_{j}}{1-q_{j}}.\label{defci+1}%
\end{equation}
The phase factor $\sigma_{j}$ reflects an interplay between the length of a
subinterval, the value of $c$ on the interval, and the wave number $\omega$.
This phase factor will be key to the analysis of the growth behaviour of the
sequence $\left(  Q_{\ell}\right)  _{\ell}$ defined later in
(\ref{defamnqmntot}) and mentioned already in (\ref{Cstabintro}).
\end{remark}

\paragraph{Variational formulation}

Let $V:=H^{1}\left(  \Omega\right)  $. Find $u\in V$ such that%
\[
\int_{\Omega}\left(  u^{\prime}\bar{v}^{\prime}-\left(  \frac{\omega}%
{c}\right)  ^{2}u\bar{v}\right)  -\operatorname*{i}\sum_{x\in\left\{
-1,1\right\}  }\frac{\omega}{c}u(x)\bar{v}(x)=
%\sum_{x\in\left\{  -1,1\right\} }g_{x}\bar{v}(x)\qquad\forall v\in V.
g_1\bar{v}(-1)+ g_2 \bar{v}(1) \qquad\forall v\in V.
\]

\subsection{Green's Function and Piecewise Constant Wave Numbers}

For $1\leq j\leq n+1$, we set $u_{j}:=\left.  u\right\vert _{\tau_{j}}$. The
solution of the homogeneous equation on the interval $\tau_{i}$ can be written
in the form%
\begin{equation}
u_{j}\left(  x\right)  =A_{j}\operatorname*{e}\nolimits^{\operatorname*{i}%
\frac{\omega}{c_{j}}x}+B_{j}\operatorname*{e}\nolimits^{-\operatorname*{i}%
\frac{\omega}{c_{j}}x}\qquad\text{on }\tau_{j}\label{localhomsol}%
\end{equation}
for $A_{j},B_{j}$, $1\leq j\leq n+1$. These coefficients are determined by the
transmission conditions that $u$ and $u^{\prime}$ are continuous at the inner
mesh points $x_{j}$, $1\leq j\leq n,$ and by the boundary conditions at
$x_{0}$ and $x_{n+1}$. The boundary conditions lead to
\begin{equation}
A_{1}=\operatorname*{i}\frac{c_{1}}{2\omega}\operatorname*{e}%
\nolimits^{\operatorname*{i}\frac{\omega}{c_{1}}}g_{1}\quad\text{and\quad
}B_{n+1}=\operatorname*{i}\frac{c_{n+1}}{2\omega}\operatorname*{e}%
\nolimits^{\operatorname*{i}\frac{\omega}{c_{n+1}}}g_{2}.\label{A1Bn}%
\end{equation}
The remaining coefficients
\begin{equation}
\mathbf{x}^{\left(  2n\right)  }:=\left(  B_{1},A_{2},B_{2},\ldots,A_{n}%
,B_{n},A_{n+1}\right)  ^{\intercal}\label{def:xu}%
\end{equation}
are the solution of a system of linear equations. Let
\begin{equation}
\mathbf{r}^{\left(  2n\right)  }=\frac{\operatorname*{i}}{2\omega}\left(
\operatorname*{e}\nolimits^{\operatorname*{i}\frac{\omega}{c_{1}}}%
g_{1},0,\ldots,0,\operatorname*{e}\nolimits^{\operatorname*{i}\frac{\omega
}{c_{n}}}g_{2}\right)  ^{\intercal}\in\mathbb{R}^{2n}.\label{defr}%
\end{equation}
We define the diagonal matrix $\mathbf{D}^{\left(  2n\right)  }\in
\mathbb{C}^{\left(  2n\right)  \times\left(  2n\right)  }$ by%
\begin{equation}
\mathbf{D}^{\left(  2n\right)  }=\operatorname*{diag}\left[  \alpha_{1,1}%
\sqrt{c_{1}},\frac{\sqrt{c_{2}}}{\alpha_{2,1}},\alpha_{2,2}\sqrt{c_{2}}%
,\frac{\sqrt{c_{3}}}{\alpha_{3,2}},\ldots,\alpha_{n,n}\sqrt{c_{n}},\frac
{\sqrt{c_{n+1}}}{\alpha_{n+1,n}}\right]  ,\label{defD}%
\end{equation}
with
\[
\alpha_{\ell,j}:=\exp\left(  \operatorname*{i}\frac{\omega}{c_{\ell}}x_{j}\right)
,\qquad1\leq j\leq n,~1\leq \ell\leq n+1,
\]
and the Green's function%
\begin{equation}
\mathbf{M}_{\mathrm{Green}}^{\left(  2n\right)  }:=\left(  \mathbf{M}^{\left(
2n\right)  }\right)  ^{-1}=%
\begin{bmatrix}
\mathbf{W}^{\left(  1\right)  } & \mathbf{N}^{\left(  1\right)  } & \mathbf{0}
& \dots & \mathbf{0}\\
\left(  \mathbf{N}^{\left(  1\right)  }\right)  ^{\intercal} & \mathbf{W}%
^{\left(  2\right)  } & \ddots & \ddots & \vdots\\
\mathbf{0} & \ddots & \ddots &  & \mathbf{0}\\
\vdots & \ddots &  &  & \mathbf{N}^{\left(  n-1\right)  }\\
\mathbf{0} & \dots & \mathbf{0} & \left(  \mathbf{N}^{\left(  n-1\right)
}\right)  ^{\intercal} & \mathbf{W}^{\left(  n\right)  }%
\end{bmatrix}^{-1} \in \mathbb{C}^{(2n)\times (2n)}.\label{defM2n-2}%
\end{equation}
The $2\times2$ submatrices $\mathbf{W}^{\left(  j\right)  }$ and
$\mathbf{N}^{\left(  j\right)  }$ are given by%
\[
\mathbf{W}^{\left(  j\right)  }:=%
\begin{bmatrix}
q_{j} & \sqrt{1-q_{j}^{2}}\\
\sqrt{1-q_{j}^{2}} & -q_{j}%
\end{bmatrix}
,\text{\quad}\mathbf{N}^{\left(  j\right)  }:=%
\begin{bmatrix}
0 & 0\\
-\frac{1}{\sqrt{\sigma_{j}}} & 0
\end{bmatrix}
\text{ with }\sqrt{\sigma_{j}}:=\exp\left(  -\operatorname*{i}\frac
{h_{j+1}\omega}{c_{j+1}}\right)  .
\]

\begin{remark}
$\mathbf{M}^{\left(  2n\right)  }$ is symmetric, but not Hermitian.
\end{remark}

The derivation of the following lemma can be found in Section \ref{ProofsRep}.

\begin{lemma}
\label{LemLGS}The remaining coefficients $\mathbf{x}^{\left(  2n \right)  }$
are given by%
\begin{equation}
\mathbf{x}^{\left(  2n \right)  }=\mathbf{D}^{\left(  2n \right)  }%
\mathbf{M}_{\mathrm{Green}}^{\left(  2n \right)  }\mathbf{D}^{\left(  2n
\right)  }\mathbf{r}^{\left(  2n \right)  }. \label{diagonalred}%
\end{equation}

\end{lemma}

\begin{remark}
The matrix $\mathbf{W}^{\left(  i\right)  }$ corresponds to a reflection; it
holds $\left(  \mathbf{W}^{\left(  i\right)  }\right)  ^{2}=\mathbf{I}$ and
the eigenvalues of $\mathbf{W}^{\left(  i\right)  } $ are $-1$ and $1$.
\end{remark}

For later use, we define the reduced matrix $\mathbf{M}^{\left(  2n-1\right)
}$ which arises by removing the last row and last column of $\mathbf{M}%
^{\left(  2n\right)  }$. Explicitly, it holds%
\[
\mathbf{M}^{\left(  2n-1\right)  }:=\left[
\begin{array}
[c]{rr}%
\mathbf{M}^{\left(  2n-2\right)  } & \mathbf{b}^{\left(  2n-2\right)  }\\
\left(  \mathbf{b}^{\left(  2n-2\right)  }\right)  ^{\intercal} & q_{n}%
\end{array}
\right]  \quad\text{with\quad}\mathbf{b}^{\left(  2n-2\right)  }:=\left(
0,0,\ldots,0,-\frac{1}{\sqrt{\sigma_{n-1}}}\right)  ^{\intercal}.
\]

One key role for the analysis of the solution operator of (\ref{Helmpwconst})
is played by the derivation of a representation of the entries of
$\mathbf{M}_{\mathrm{Green}}^{\left(  2n\right)  }$. We denote by
$\mathbf{M}^{\left(  2n,i,j\right)  }$ the matrix which arises if the $i$-th
row and the $j$-th column of $\mathbf{M}^{\left(  2n\right)  }$ are removed.
According to Cramer's rule we have%
\[
\left(  \mathbf{M}_{\mathrm{Green}}^{\left(  2n\right)  }\right)
_{i,j}=\left(  -1\right)  ^{i+j}\frac{\det\mathbf{M}^{\left(  2n,i,j\right)
}}{\det\mathbf{M}^{\left(  2n\right)  }}.
\]
To express the determinant of $\mathbf{M}^{\left(  2n\right)  }$ as a product
of terms with positive modulus we introduce some notation. We define the sequence $Q_{m}$ by the
recursion%
%TCIMACRO{\TeXButton{defamnqmntot}{\begin{subequations}
%\label{defamnqmntot}
%\end{subequations}}}%
%BeginExpansion
\begin{subequations}
\label{defamnqmntot}
\end{subequations}%
%EndExpansion%
\begin{equation}%
\begin{array}
[c]{l}%
Q_{1}=\frac{q_{1}}{\sigma_{1}},\\
Q_{j}=\dfrac{q_{j}+Q_{j-1}}{\sigma_{j}\left(  1+q_{j}Q_{j-1}\right)  }%
,\quad2\leq j\leq n.
\end{array}
\tag{%
%TCIMACRO{\TeXButton{defamnqmntot}{\ref{defamnqmntot}}}%
%BeginExpansion
\ref{defamnqmntot}%
%EndExpansion
a}\label{defamnqmna}%
\end{equation}
For later use, we define the quantity $Q_{j}^{\prime}$ of same modulus by
\begin{equation}
Q_{j}^{\prime}:=\sigma_{j}Q_{j}\tag{%
%TCIMACRO{\TeXButton{defamnqmntot}{\ref{defamnqmntot}}}%
%BeginExpansion
\ref{defamnqmntot}%
%EndExpansion
b}\label{defamnqmnb}%
\end{equation}
Formally, we set $Q_{0}=0$. Note that
$Q_{j}^{\prime}$ is independent of $\sigma_{j}$ and $Q_{j}$ depends on
$\left(  \sigma_{i}\right)  _{i=1}^{j}$, $\left(  q_{i}\right)  _{i=1}^{j}$.
This allows to define%
\begin{equation}
\tilde{p}_{n}\left(  \boldsymbol{\sigma},\mathbf{q}\right)  :=\prod
\nolimits_{j=1}^{n-1}\left(  1+q_{j+1}Q_{j}\right)  \label{defpntilde}%
\end{equation}
Here, $\boldsymbol{\sigma}=\left(  \sigma_{j}\right)  _{j=1}^{n-1}$ and
$\mathbf{q}=\left(  q_{j}\right)  _{j=1}^{n}$ in $\tilde{p}_{n}$. The proof of
the following lemma is postponed to Section \ref{ProofsRep}.

\begin{lemma}
\label{LemDet}For $n\geq2$, it holds%
\begin{align*}
\det\mathbf{M}^{\left(  2n\right)  }  &  =\left(  -1\right)  ^{n}\tilde{p}%
_{n},\\
\det\mathbf{M}^{\left(  2n-1\right)  }  &  =-\sigma_{n}Q_{n}\det
\mathbf{M}^{\left(  2n\right)  }.
\end{align*}

\end{lemma}

In view of the special structure of the right-hand side \eqref{defr}, we are
particularly interested in the first and last column of $\mathbf{M}%
_{\mathrm{Green}}^{\left(  2n\right)  }$, denoted by $\left(  \mathbf{M}%
_{\mathrm{Green}}^{\left(  2n\right)  }\right)  _{\ast,1}$, $\left(
\mathbf{M}_{\mathrm{Green}}^{\left(  2n\right)  }\right)  _{\ast,2n}$. We
investigate the last column, i.e., $j=2n$. By symmetry, the computations for
the first column are equivalent. From \eqref{formulacof} we obtain
\begin{equation}
\det\mathbf{M}^{\left(  2n,i,2n\right)  }=\left(  {\prod_{j=\left\lfloor
\frac{i+2}{2}\right\rfloor }^{n}}\sqrt{1-q_{j}^{2}}\right)  \left(
\prod_{j=\left\lfloor \frac{i+1}{2}\right\rfloor }^{n-1}\left(  -\frac
{1}{\sqrt{\sigma_{j}}}\right)  \right)  \det\mathbf{M}^{\left(  i-1\right)
}.\label{detrep1}%
\end{equation}
The following corollary is a direct consequence of Lemma \ref{LemDet} and \eqref{detrep1}. It allows us to write the entries of the Green's function $\mathbf{M}%
_{\mathrm{Green}}$ in terms of the sequence $(Q_{j})$, leading to a stability
estimate depending on $(Q_{j})$ (Theorem \ref{Theo:StabEst}).

\begin{corollary}
\label{CorrepR}The entries of the last column of $\mathbf{M}_{\mathrm{Green}%
}^{\left(  2n\right)  }$ can be written in the form%
\begin{equation}
\left(  \mathbf{M}_{\mathrm{Green}}^{\left(  2n\right)  }\right)
_{i,2n}=\left(  -1\right)  ^{i+1}\left(  {\prod_{\ell=\left\lfloor \frac
{i+2}{2}\right\rfloor }^{n}}\frac{\sqrt{1-q_{\ell}^{2}}}{\left(  1+q_{\ell
}Q_{\ell-1}\right)  \sqrt{\sigma_{\ell-1}}}\right)  \times%
\begin{cases}
\sigma_{\frac{i}{2}}Q_{\frac{i}{2}} & i\text{ is even,}\\
\sqrt{\sigma_{\frac{i-1}{2}}} & i\text{ is odd.}%
\end{cases}
\label{repsol}%
\end{equation}
For the modulus of the entries of the odd rows it holds%
\begin{equation}
\left\vert \left(  \mathbf{M}_{\mathrm{Green}}^{\left(  2n\right)  }\right)
_{2n-2m+1,2n}\right\vert =\left\vert {\prod_{\ell=n-m+1}^{n}}\frac
{\sqrt{1-q_{\ell}^{2}}}{\left(  1+q_{\ell}Q_{\ell-1}\right)  }\right\vert
,\qquad1\leq m\leq n,\label{Estoddrows}%
\end{equation}
and for the even rows:%
\begin{align*}
\left\vert \left(  \mathbf{M}_{\mathrm{Green}}^{\left(  2n\right)  }\right)
_{2n,2n}\right\vert  &  =\left\vert Q_{n}\right\vert ,\\
\left\vert \left(  \mathbf{M}_{\mathrm{Green}}^{\left(  2n\right)  }\right)
_{2n-2m,2n}\right\vert  &  =\left\vert \left(  \mathbf{M}_{\mathrm{Green}%
}^{\left(  2n\right)  }\right)  _{2n-2m+1,2n}\right\vert \left\vert
Q_{n-m}\right\vert ,\qquad1\leq m\leq n-1.
\end{align*}

\end{corollary}

The estimates on the coefficients in $\mathbf{M}_{\mathrm{Green}}^{\left(
2n\right)  }$ is based on the representation (\ref{Estoddrows}), i.e., we have
to estimate the expression%
\begin{equation}
G_{n,m}:=\left\vert {\prod_{\ell=n-m+1}^{n}}\frac{\sqrt{1-q_{\ell}^{2}}%
}{\left(  1+q_{\ell}Q_{\ell-1}\right)  }\right\vert ,\qquad1\leq m\leq
n,\label{defgnm}%
\end{equation}
with $Q_{j}$ as in (\ref{defamnqmna}).

\begin{lemma}
\label{LemGnm}

\begin{enumerate}
\item[(i)] Let $q_{j}\in\left[  -q,q\right]  $ for all $1\leq j\leq n$ and some
$0<q<1$. Then%
\[
G_{n,m}\leq\frac{1}{\sqrt{1-\left\vert Q_{n-m}\right\vert ^{2}}},\qquad1\leq
m\leq n.
\]

\item[(ii)] If additionally $Q_{n-1}=-q_{n}$ holds, we have
\begin{equation}
G_{n,m}=\frac{1}{\sqrt{1-|Q_{n-m}|^{2}}},\qquad1\leq m\leq n.
\label{choiceofsq}%
\end{equation}

\end{enumerate}
\end{lemma}

\proof

\begin{enumerate}
\item[(i)] We prove the statement by induction on $m$. For $m=1$, we employ
simple calculus and obtain%
\[
G_{n,1}\leq\max_{q_{n}\in\left[  -q,q\right]  }\left\vert \frac{\sqrt
{1-q_{n}^{2}}}{1+q_{n}Q_{n-1}}\right\vert =\max_{x\in\left[  0,q\right]
}\frac{\sqrt{1-x^{2}}}{1-x\left\vert Q_{n-1}\right\vert }\leq\frac{1}%
{\sqrt{1-\left\vert Q_{n-1}\right\vert ^{2}}}.
\]
Next we assume that the estimate holds for $G_{n,\ell}$, $1\leq\ell\leq m-1$.
For $m$, we get%
\begin{equation}%
\begin{split}
G_{n,m} &  =\frac{\sqrt{1-q_{n-m+1}^{2}}}{\left\vert 1+q_{n-m+1}%
Q_{n-m}\right\vert }{\prod\nolimits_{\ell=n+2-m}^{n}}\frac{\sqrt{1-q_{\ell
}^{2}}}{\left\vert 1+q_{\ell}Q_{\ell-1}\right\vert }\\
&  \leq\frac{\sqrt{1-q_{n-m+1}^{2}}}{\left\vert 1+q_{n-m+1}Q_{n-m}\right\vert
}\frac{1}{\sqrt{1-\left\vert Q_{n-m+1}\right\vert ^{2}}}\\
&  =\frac{\sqrt{1-q_{n-m+1}^{2}}}{\sqrt{\left\vert 1+q_{n-m+1}Q_{n-m}%
\right\vert ^{2}-\left\vert q_{n-m+1}+Q_{n-m}\right\vert ^{2}}}\\
&  =\frac{1}{\sqrt{1-\left\vert Q_{n-m}\right\vert ^{2}}}.
\end{split}
\label{GnmInduction}%
\end{equation}

\item[(ii)] Since $Q_{n-1}=-q_{n}$, we can compute $G_{n,1}$ directly
\[
G_{n,1}=\left\vert \frac{\sqrt{1-q_{n}^{2}}}{1+q_{n}Q_{n-1}}\right\vert
%=\left\vert \frac{\sqrt{1-q_{n}^{2}}}{1+q_{n}\frac{Q_{n-1}^{\prime}}{\sigma_{n-1}}}\right\vert\overset{\text{(\ref{setsigmaq})}}{=}
=\frac{1}{\sqrt{1-|Q_{n-1}|^{2}}}.
\]
The induction step follows analogously to (\ref{GnmInduction}).
\end{enumerate}

\endproof

With the estimate of the entries of $\mathbf{M}_{\mathrm{Green}}^{\left(
2n\right)  } $ of Lemma \ref{LemGnm}, we can state the following stability result.

\begin{theorem}
\label{Theo:StabEst}The solution $u$ of (\ref{Helmpwconst}) with piecewise
constant wave speed $c$ as in \eqref{defc} satisfies%
\begin{equation}
\left\Vert u^{\left(  k\right)  }\right\Vert _{L^{2}\left(  \Omega\right)
}\leq4\frac{c_{\max}}{c_{\min}^{k}}\omega^{k-1}\max\left\{  \left\vert
g_{1}\right\vert ,\left\vert g_{2}\right\vert \right\}  \max_{1\leq j\leq
n}\frac{1}{\sqrt{1-\left\vert Q_{j}\right\vert ^{2}}} \label{estuk}%
\end{equation}
for $k=0,1,2$.
\end{theorem}

\begin{remark}
If $c$ is constant, the qualitative frequency dependence in the stability
estimate \eqref{estuk} coincides with known results (cf. \cite{MelenkDiss,
CummingsFeng2006}).
\end{remark}

\begin{remark}
The final result (Theorem \ref{Th:LastTh}) is proved by a further analysis of
the sequence $(Q_{j})$ or, more precisely, of the distance $1-\left\vert
Q_{j}\right\vert ^{2}$. We emphasize that the magnitude of the last term in
\eqref{estuk} does not necessarily increase in $\omega$, but depends on the
interplay of $h_{j},c_{j}$ and $\omega$ via the phase factor $\sigma_{j}$
defined in \eqref{defalphaij}(resonance effect) .
\end{remark}

\proof From (\ref{localhomsol}) we obtain for any $k=0,1,2$%
\[
\left\Vert u^{\left(  k\right)  }\right\Vert _{L^{2}\left(  \Omega\right)
}\leq\sqrt{2}\max_{1\leq j\leq n+1}\left(  \frac{\omega}{c_{j}}\right)
^{k}\left(  \left\vert A_{j}\right\vert +\left\vert B_{j}\right\vert \right)
.
\]
The coefficients $A_{j}$, $B_{j}$ are contained in the solution vector
$\mathbf{x}^{\left(  2n\right)  }$ and can be written in the form (cf.
(\ref{diagonalred}))%
\begin{align*}
\mathbf{x}_{j}^{\left(  2n\right)  } &  =\left(  \left(  \mathbf{D}^{\left(
2n\right)  }\right)  \mathbf{M}_{\mathrm{Green}}^{\left(  2n\right)  }\left(
\mathbf{D}^{\left(  2n\right)  }\right)  \mathbf{r}^{\left(  2n\right)
}\right)  _{j}\\
&  =\frac{\operatorname*{i}}{2\omega}\left(  \alpha_{1,1}\sqrt{c_{1}%
}\operatorname*{e}\nolimits^{\operatorname*{i}\frac{\omega}{c_{1}}}%
g_{1}\mathbf{d}_{j}^{\left(  2n\right)  }\left(  \mathbf{M}_{\mathrm{Green}%
}^{\left(  2n\right)  }\right)  _{j,1}+\frac{\sqrt{c_{n+1}}}{\alpha_{n+1,n}%
}\operatorname*{e}\nolimits^{\operatorname*{i}\frac{\omega}{c_{n+1}}}%
g_{2}\mathbf{d}_{j}^{\left(  2n\right)  }\left(  \mathbf{M}_{\mathrm{Green}%
}^{\left(  2n\right)  }\right)  _{j,2n}\right)
\end{align*}
with (cf. (\ref{defD}))%
\[
\mathbf{d}^{\left(  2n\right)  }:=\left(  \alpha_{1,1}\sqrt{c_{1}},\frac
{\sqrt{c_{2}}}{\alpha_{2,1}},\alpha_{2,2}\sqrt{c_{2}},\frac{\sqrt{c_{3}}%
}{\alpha_{3,2}},\ldots,\alpha_{n,n}\sqrt{c_{n}},\frac{\sqrt{c_{n+1}}}%
{\alpha_{n+1,n}}\right)  ^{\intercal}.
\]
Hence,%
\begin{equation}
\max_{1\leq j\leq n}\max\left\{  \left\vert A_{j}\right\vert ,\left\vert
B_{j}\right\vert \right\}  \leq\frac{c_{\max}}{\omega}\max\left\{  \left\vert
g_{1}\right\vert ,\left\vert g_{2}\right\vert \right\}  \max_{1\leq i\leq
2n}\max_{j\in\left\{  1,2n\right\}  }\left\vert \left(  \mathbf{M}%
_{\mathrm{Green}}^{\left(  2n\right)  }\right)  _{i,j}\right\vert
.\label{eq:stabilityfzero}%
\end{equation}
\endproof

\begin{lemma}
\label{lem:lowerbound} A lower bound for the norm of $u^{\left(  k\right)  }$,
$k=0,1,2$, is given by%
\[
\left\Vert u^{\left(  k\right)  }\right\Vert _{L^{2}\left(  \Omega\right)
}\geq\sqrt{\frac{2}{15}h_{j}}\left(  \frac{\omega}{c_{j}}\right)  ^{k}%
\frac{\frac{\omega}{c_{j}}h_{j}}{1+\frac{\omega}{c_{j}}h_{j}}\max\left\{
\left\vert A_{j}\right\vert ,\left\vert B_{j}\right\vert \right\}  ,
\]
where $u$ and $A_{j},B_{j}$ are related through \eqref{localhomsol} and
(\ref{def:xu}).
\end{lemma}

%

%TCIMACRO{\TeXButton{Proof}{\proof}}%
%BeginExpansion
\proof
%EndExpansion
It holds%
\[
\left\Vert u^{\left(  k\right)  }\right\Vert^2_{L^{2}\left(  \Omega\right)
}\geq\left\Vert u^{\left(  k\right)  }\right\Vert^2_{L^{2}\left(  \tau
_{j}\right)  }=\left\langle \binom{A_{j}}{B_{j}},E^{\left(  k\right)  }%
\binom{A_{j}}{B_{j}}\right\rangle
\]
with the Euclidean scalar product $\left\langle \cdot,\cdot\right\rangle $,
the Hermitian matrix
\[
E_{n,m}^{\left(  k\right)  }:=\left(  \frac{\omega}{c_{j}}\right)  ^{2k}%
h_{j}\left[
\begin{array}
[c]{cc}%
1 & \left(  -1\right)  ^{k}\operatorname*{e}\nolimits^{\operatorname*{i}%
\frac{\omega}{c_{j}}\left(  2x_{j-1}+h_j\right)  }\operatorname{sinc}\left(
\frac{\omega}{c_{j}}h_{j}\right)  \\
\left(  -1\right)  ^{k}\operatorname*{e}\nolimits^{-\operatorname*{i}%
\frac{\omega}{c_{j}}\left(  2x_{j-1}+h_j\right)  }\operatorname{sinc}\left(
\frac{\omega}{c_{j}}h_{j}\right)   & 1
\end{array}
\right]  ,
\]
and the convention for the sinc function $\operatorname{sinc}\left(  x\right)
:=\left(  \sin x\right)  /x$ for $x\neq0$ and $\operatorname{sinc}\left(
0\right)  :=1$. The eigenvalues of $E^{\left(  k\right)  }$ are given by%
\[
\lambda_{1,2}^{\left(  k\right)  }=\left(  \frac{\omega}{c_{j}}\right)
^{2k}h_{j}\left(  1\pm\operatorname{sinc}\left(  \frac{\omega}{c_{j}}%
h_{j}\right)  \right)  .
\]
Simple calculus leads to%
\[
\left\vert \lambda_{1,2}^{\left(  k\right)  }\right\vert \geq\left(
\frac{\omega}{c_{j}}\right)  ^{2k}h_{j}\times\left\{
\begin{array}
[c]{cc}%
\frac{2}{15}\left(  \frac{\omega}{c_{j}}h_{j}\right)  ^{2} &  \frac{\omega}{c_{j}}h_{j}\leq2\\
\frac{1}{2} &  \frac{\omega}{c_{j}}h_{j}\geq2
\end{array}
\right\}  \geq\frac{2}{15}h_{j}\left(  \frac{\omega}{c_{j}}\right)  ^{2k}%
\frac{\left(  \frac{\omega}{c_{j}}h_{j}\right)  ^{2}}{1+\left(  \frac{\omega
}{c_{j}}h_{j}\right)  ^{2}}%
\]
and some straightforward manipulations result in the asserted estimate.%
%TCIMACRO{\TeXButton{End Proof}{\endproof}}%
%BeginExpansion
\endproof
%EndExpansion

\section{Construction of Configurations with \textquotedblleft
good\textquotedblright\ and \textquotedblleft bad\textquotedblright\ Stability
Properties}

\label{Sec:ex} 
\begin{figure}[ptb]
\begin{subfigure}{.47\textwidth}
%\centering
\includegraphics[scale=0.47]{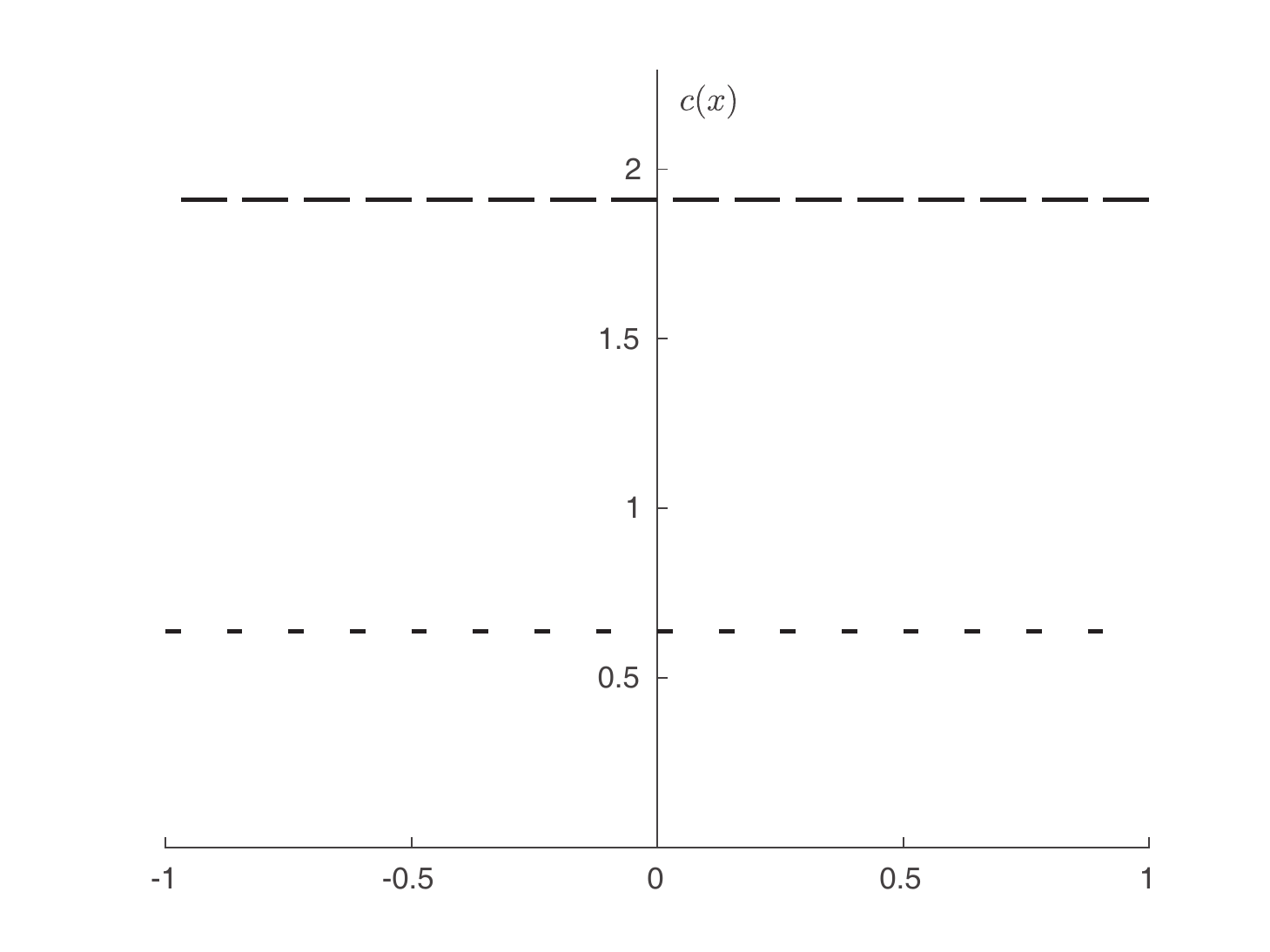}
\includegraphics[scale=0.47]{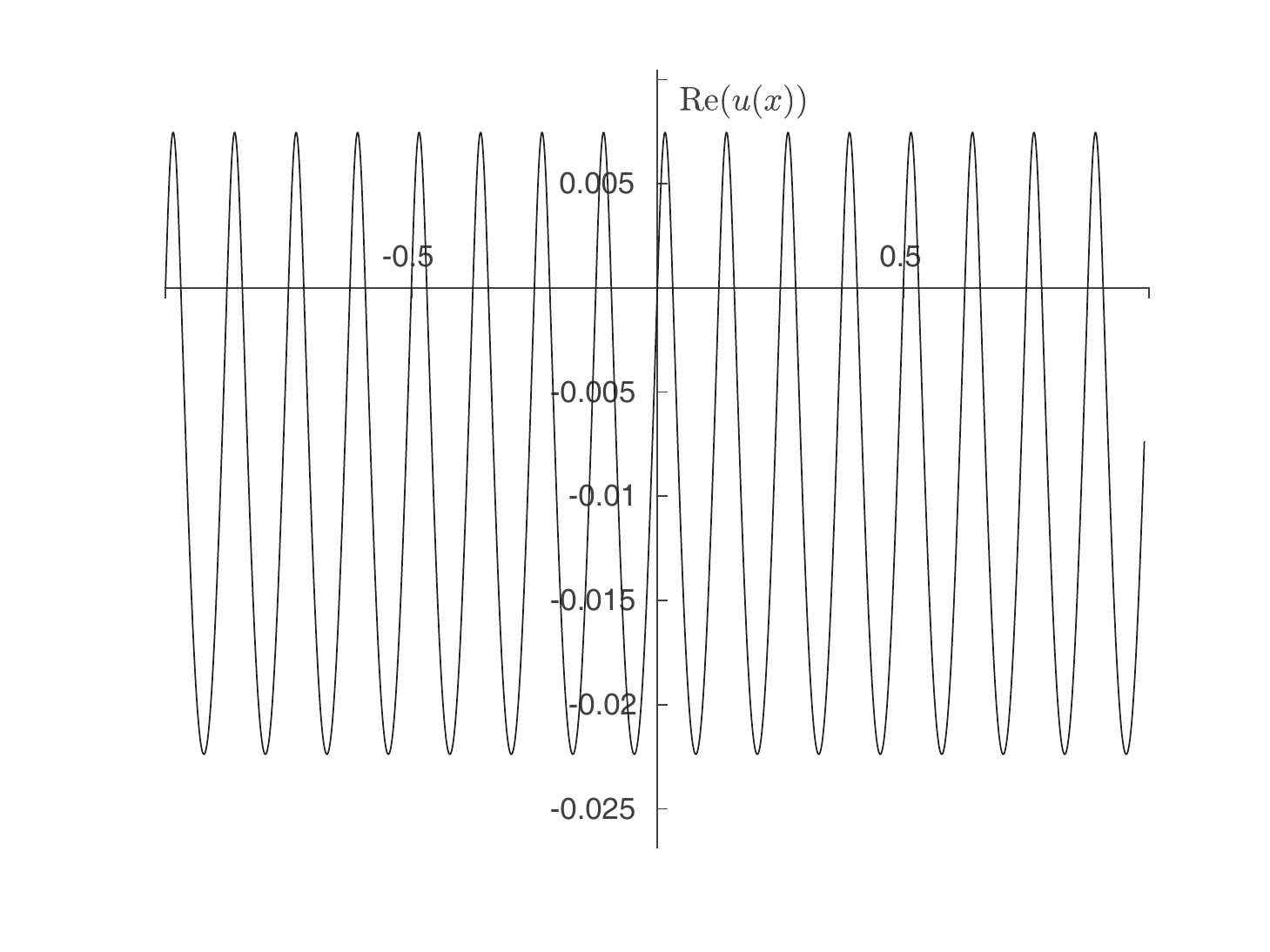}
\includegraphics[scale=0.47]{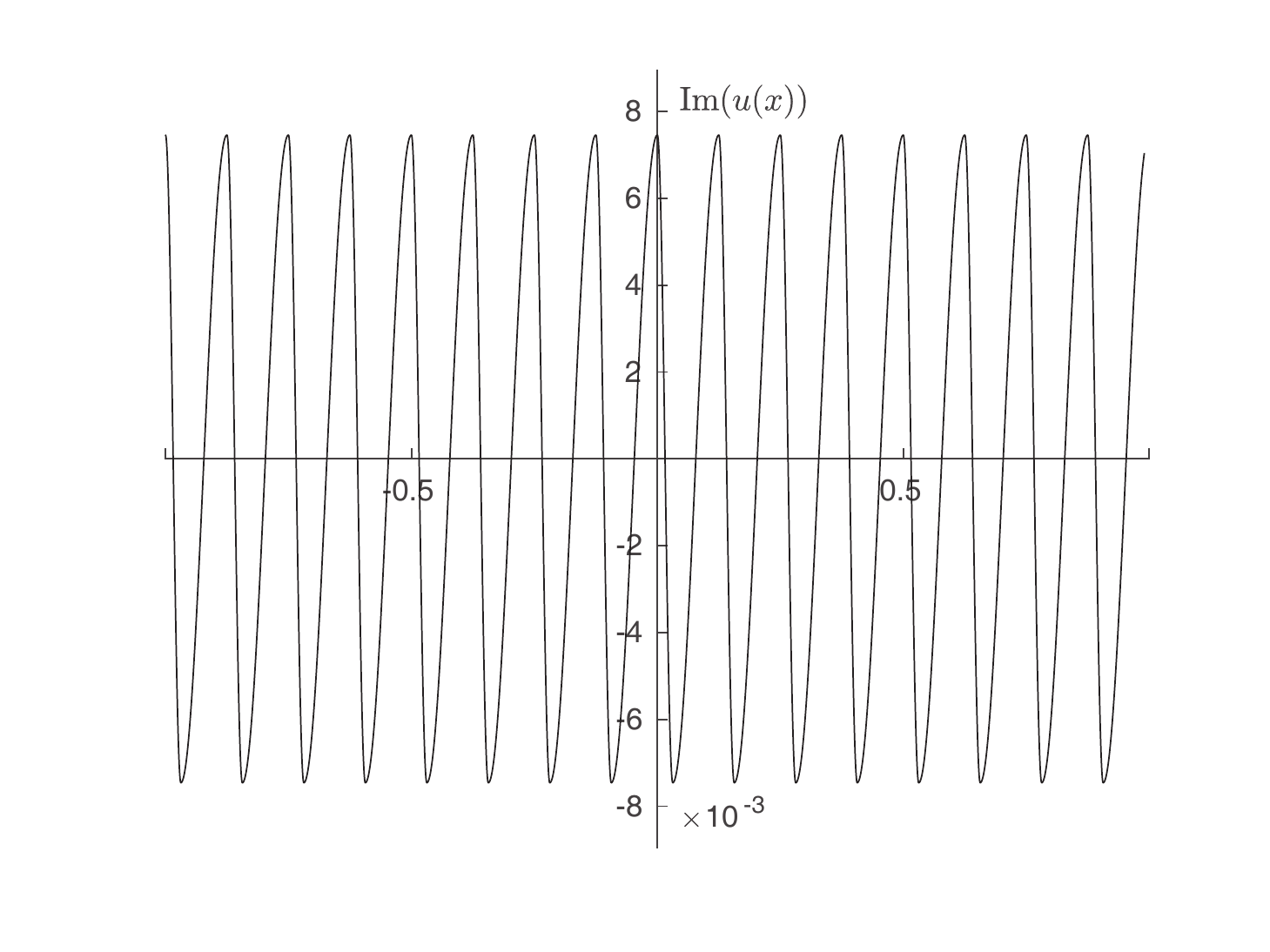}
\caption{Solution \(u\), with \(n = 32\). The stability constant is bounded from above independently of \(\omega\)}
\end{subfigure}
\hfill\begin{subfigure}{.47\textwidth}
\centering%
\includegraphics[scale=0.47]{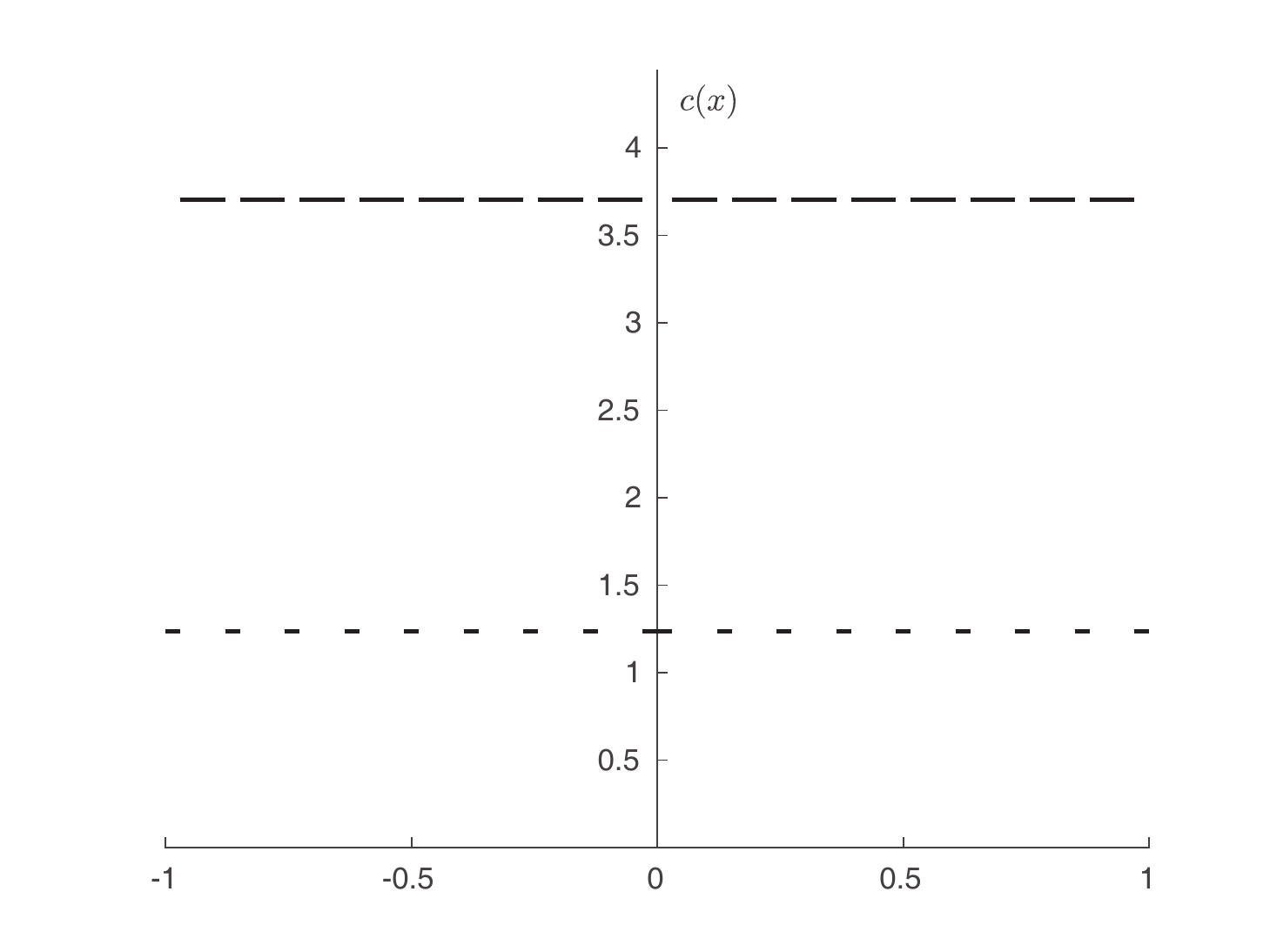}
\includegraphics[scale=0.47]{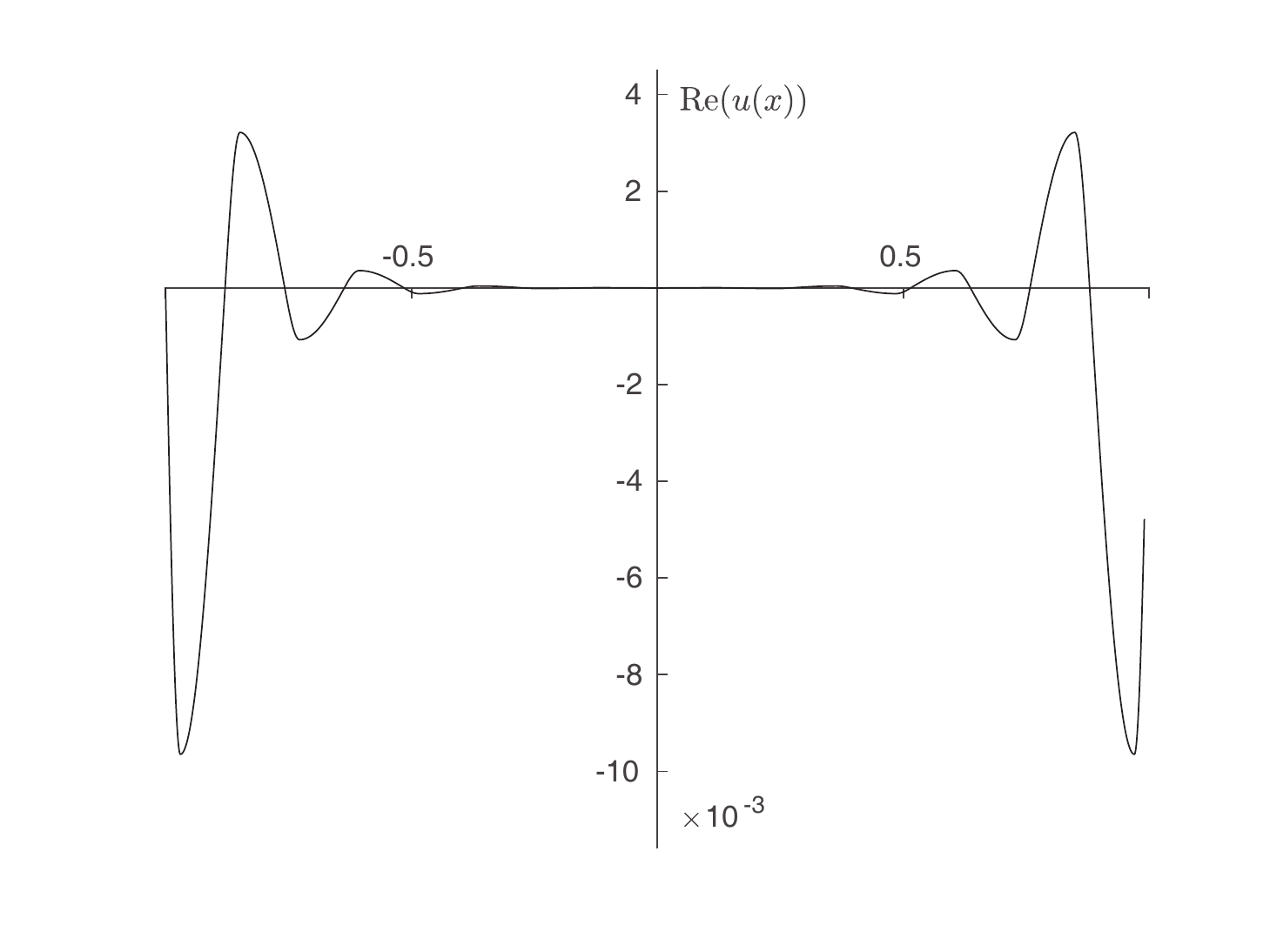}
\includegraphics[scale=0.47]{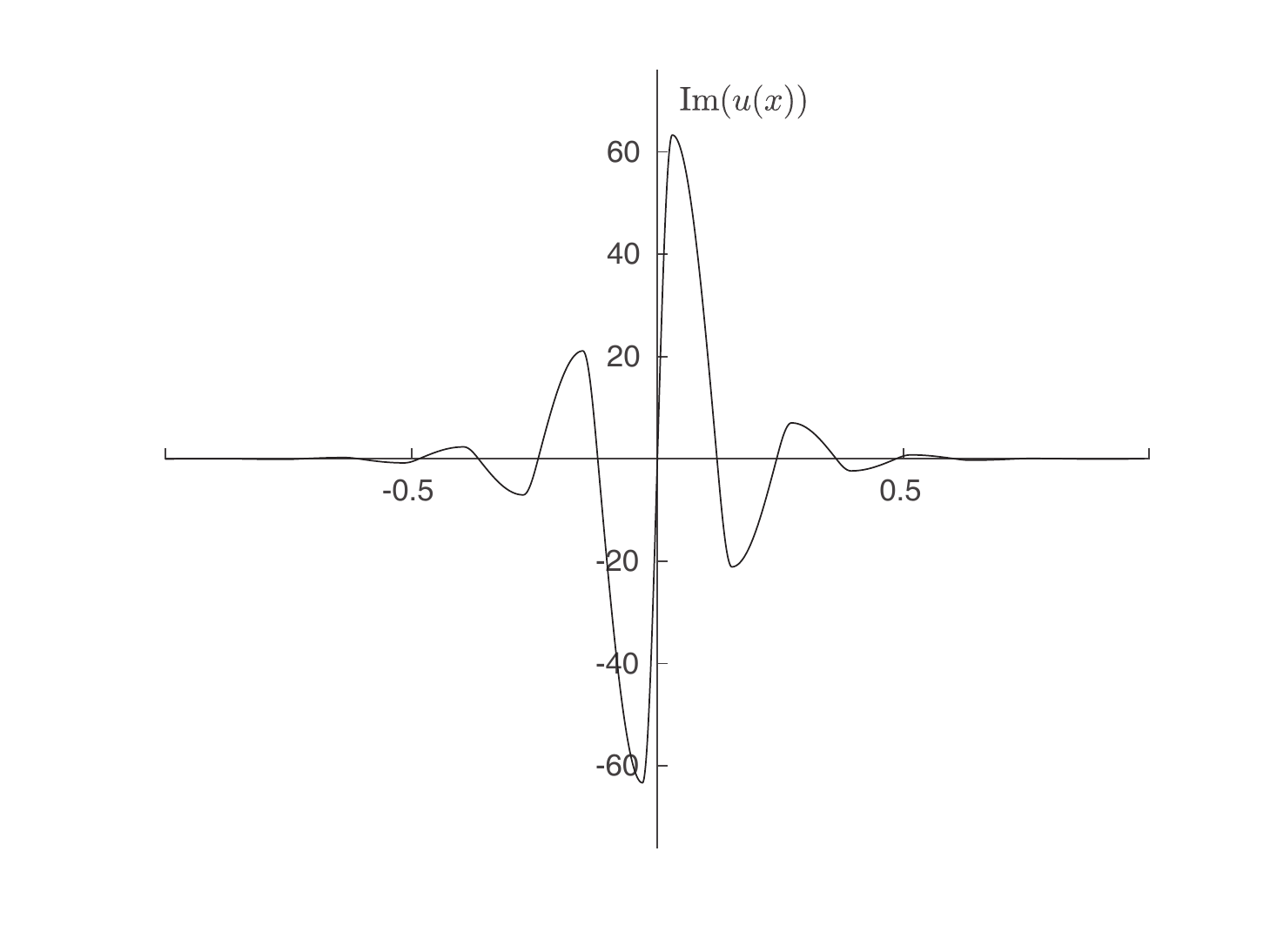}
\caption{Solution \(u\), with \(n = 33\). The stability constant grows exponentially in \(\omega\). }
\end{subfigure}
\caption{Examples for a piecewise constant wave speed and corresponding
solution of the Helmholtz problem. In both examples we set $\omega=64$,
$q=0.5$, $g_{1}=0$, $g_{2}=1$, $h_{i}$, $c$ as in \S \ref{example_osc}
and \S \ref{example_exp}, respectively.}%
\label{Fig:Examples}%
\end{figure}

We construct two slightly different configurations of $c_{j}$
and $h_{j}$, to find two very different behaviour of the growth of $|Q_{j}|$
leading to qualitatively different maximal entries of $\mathbf{M}%
_{\mathrm{Green}}^{(2n)}$. Figure \ref{Fig:Examples} shows the corresponding
solution for a specific example in both cases. In both cases, we first choose
the \textit{mathematical parameters} $\sigma_{j}$, $q_{j}$ such that the
sequence $\left(  Q_{\ell}\right)  _{\ell}$ exhibits the desired behaviour and
then determine the corresponding physical parameters.

\subsection{Well-behaved Case}

\label{example_osc}
Given an arbitrary frequency $\omega$, we construct an example of a highly
oscillating configuration, where $|Q_{j}|$ can be bounded away from 1
independently of $\omega$. Recalling the definition of $Q_{j}$, we write
$Q_{2}$ as
\[
Q_{2}=\frac{q_{2}+\frac{q_{1}}{\sigma_{1}}}{\sigma_{2}\left(  1+\frac
{q_{1}q_{2}}{\sigma_{1}}\right)  }.
\]
We see that the modulus of $Q_{2}$ is minimized, if $q_{1},q_{2}$ and
$\sigma_{1}$ are chosen such that $q_{2}+\frac{q_{1}}{\sigma_{1}}=0$. This is
true, for example, if $q_{1}=-q_{2}$ and $\sigma_{1}=1$. If indeed $Q_{2}=0$,
the same choice for $q_{3},q_{4}$ and $\sigma_{3}$ leads to $Q_{4}=0$. As a
result of this observation, we define the oscillating, piecewise constant wave
speed by%
\begin{equation}
c_{j}=%
\begin{cases}
c\left(  1-q\right)   & j\text{ odd,}\\
c\left(  1+q\right)   & j\text{ even,}%
\end{cases}
\qquad1\leq j\leq n+1,\label{defwavespeedex}%
\end{equation}
for some $c>0$ and $q\in\left(  0,1\right)  $, so that $q=\frac{c_{\max
}-c_{\min}}{c_{\max}+c_{\min}}$. Indeed, we then obtain $q_{j}=\left(
-1\right)  ^{j+1}q$. In order to have $\sigma_{j}=1$ for all $j$ (cf.
\eqref{defalphaij}), we consider the case $\frac{\omega h_{j}}{c_{j}}\gtrsim1$
and step widths $h_{j}$ which satisfy%
\begin{equation}
\frac{\omega h_{j}}{c_{j}\pi}\in\mathbb{N}\text{ for all }%
j.\label{condgoodcase}%
\end{equation}
To be concrete, we choose $n$ odd and
\[
h_{j}:=c_{j}\frac{\pi}{\omega}%
\]
The side condition
\[
\sum_{j=1}^{n+1}h_{j}=\sum_{j=1}^{n+1}c_{j}\frac{\pi}{\omega}=2
\]
leads to
\[
c=\frac{2}{n+1}\frac{\omega}{\pi},\text{\quad}h_{j}=\frac{2}{n+1}\times%
\begin{cases}
1-q & j\text{ odd,}\\
1+q & j\text{ even,}%
\end{cases}
%\qquad 1\leq j\leq n+1,
\]
and
\[
x_{j}=%
\begin{cases}
-1+\dfrac{2\left(  j-q\right)  }{n+1} & j\text{ odd,}\\
-1+\dfrac{2j}{n+1} & j\text{ even,}%
\end{cases}
\qquad1\leq j\leq n+1.
\]
In any case, condition (\ref{condgoodcase}) implies%
\[
Q_{j}=%
\begin{cases}
q & j\text{ odd,}\\
0 & j\text{ even,}%
\end{cases}
\]
and%
\[
\frac{\sqrt{1-q_{j}^{2}}}{1+q_{j}Q_{j-1}}=%
\begin{cases}
\frac{1}{\sqrt{1-q^{2}}} & j\text{ even,}\\
\sqrt{1-q^{2}} & j\text{ odd.}%
\end{cases}
\]
Hence, the product of two subsequent factors in $G_{n,m}$ equals $1$, so that%
\[
G_{n,m}\leq\frac{1}{\sqrt{1-q^{2}}}.
\]
Note that%
\[
\frac{1}{\sqrt{1-q^{2}}}=\frac{1}{2}\left(  \sqrt{\frac{c_{\max}}{c_{\min}}%
}+\sqrt{\frac{c_{\min}}{c_{\max}}}\right)  \leq\sqrt{\frac{c_{\max}}{c_{\min}%
}}%
\]
From Corollary \ref{CorrepR}, we conclude that%
\[
\left\vert \left(  \mathbf{M}_{\mathrm{Green}}^{\left(  2n\right)  }\right)
_{j,2n}\right\vert \leq\sqrt{\frac{c_{\max}}{c_{\min}}}.
\]
By employing the same arguments to the first column of $\mathbf{M}%
_{\mathrm{Green}}^{\left(  2n\right)  }$, we get $\left\vert \left(
\mathbf{M}_{\mathrm{Green}}^{\left(  2n\right)  }\right)  _{j,1}\right\vert
\leq\sqrt{\frac{c_{\max}}{c_{\min}}}$. Hence, Theorem \ref{Theo:StabEst}
implies%
\[
\left\Vert u^{\left(  k\right)  }\right\Vert _{L^{2}\left(  \Omega\right)
}\leq\frac{c_{\max}^{3/2}}{c_{\min}^{k+1/2}}\omega^{k-1}\max\left\{
\left\vert g_{1}\right\vert ,\left\vert g_{2}\right\vert \right\}  .
\]
This shows that the known bounds in \cite{ChaumontFreletDiss} and
\cite{GrahamSauter2017}, which grow exponentially in the number of jumps, are
very pessimistic for this example. We summarise the findings of this example.
Let $\omega>0$ be given and let $n=O(\omega)$ be odd. Choose some $q\in\left[
0,1\right[  $. Define%
\[
c:=\frac{2}{n+1}\frac{\omega}{\pi}=O(1),\text{\quad}h_{j}=\frac{2}{n+1}\times%
\begin{cases}
1-q & j\text{ odd,}\\
1+q & j\text{ even,}%
\end{cases}
,
\]
and the piecewise constant wave speed by
\[
c_{j}=\left\{
\begin{array}
[c]{ll}%
c\left(  1-q\right)   & j\text{ odd,}\\
c\left(  1+q\right)   & j\text{ even.}%
\end{array}
\right.
\]
The stability constant for problem (\ref{Helmpwconst}) is bounded independent
of $\omega$ as in the case of constant wave speed.

\subsection{Critical Case}

\label{example_exp}
%\begin{example}
For a given frequency $\omega$, we construct a class of configurations of
$c_{j}$ and $h_{j}$ which are in resonance, i.e. where we observe a stability
constant which grows exponentially in the frequency. However, we remark that
these configurations seem to be very rare, since they depend on a
specific relation between $c_{j}$, $h_{j}$ and $\omega$ in a very sensitive way. Recall the notation
as in (\ref{defamnqmntot}) and $G_{n,m}$ is defined as in (\ref{defgnm}). If
$Q_{n-1}=-q_{n}$ we know from Lemma \ref{LemGnm} b) that
\begin{equation}
G_{n,\ell}=\frac{1}{\sqrt{1-|Q_{n-\ell}|^{2}}},\qquad1\leq\ell\leq
n\text{.}\label{EqGexact}%
\end{equation}
We investigate the question whether some $|Q_{k}|^{2}$ can become
exponentially close to $1$ with respect to growing $k$. Let $n=2k$ be even. We
will choose the first $k-1$ entries of $(\sigma_{j})_{j=1}^{n}$ such that
$|Q_{j}|$ increases, and the last entries will be chosen such that $|Q_{j}|$
decreases until $Q_{n-1}=-q_{n}$. We employ the same ansatz for the perfectly
oscillating wave speed $c$ as in (\ref{defwavespeedex}), so that
$q=\frac{c_{\max}-c_{\min}}{c_{\max}+c_{\min}}\in\left(  0,1\right)  $ and
$q_{j}=\left(  -1\right)  ^{j+1}q$. We choose the phase factors $\sigma
_{j}=\exp\left(  -2\operatorname*{i}\frac{\omega h_{j+1}}{c_{j+1}}\right)  $,
and will adjust the value of the constant $c$ and $(h_{j})_{j=1}^{n+1}$. We
choose the first $k-1$ phase factors $\sigma_{j}$ according to
\[
\sigma_{j}:=\operatorname*{sign}\left(  q_{j}q_{j+1}\right)  =-1,\qquad
\forall1\leq j\leq k-1.
\]
In Lemma \ref{estA1Q1}, we will prove that then
\[
Q_{j}=(-1)^{j}\frac{\left(  1+q\right)  ^{j}-\left(  1-q\right)  ^{j}}{\left(
1+q\right)  ^{j}+\left(  1-q\right)  ^{j}},\qquad1\leq j\leq k-1.
\]
and
\[
Q_{k}=\frac{(-1)^{k-1}}{\sigma_{k}}\frac{\left(  1+q\right)  ^{k}-\left(
1-q\right)  ^{k}}{\left(  1+q\right)  ^{k}+\left(  1-q\right)  ^{k}}.
\]
Now, setting $\sigma_{k}=1$ yields
\begin{align*}
Q_{k+1} &  =\frac{1}{\sigma_{k+1}}\dfrac{q_{k+1}+Q_{k}}{1+q_{k+1}Q_{k}}%
=\frac{(-1)^{k-1}}{\sigma_{k+1}}\dfrac{-q+|Q_{k}|}{1-q|Q_{k}|}\\
&  =\frac{(-1)^{k-1}}{\sigma_{k+1}}\frac{\left(  1+q\right)  ^{k-1}-\left(
1-q\right)  ^{k-1}}{\left(  1+q\right)  ^{k-1}+\left(  1-q\right)  ^{k-1}}.
\end{align*}
With the choice $\sigma_{j}=-1,$ for all $k+1\leq j\leq n,$ and proceeding in
an analogue way we obtain for $1\leq j\leq k$%
\[
Q_{k+j}=(-1)^{k-j+1}\frac{\left(  1+q\right)  ^{k-j}-\left(  1-q\right)
^{k-j}}{\left(  1+q\right)  ^{k-j}+\left(  1-q\right)  ^{k-j}}.
\]
This leads finally to (note that $q_{n}=-q$ for even $n$)
\[
Q_{n-1}=q.
\]
Hence,%
\[
\left\vert Q_{k}\right\vert =\frac{\left(  1+q\right)  ^{k}-\left(
1-q\right)  ^{k}}{\left(  1+q\right)  ^{k}+\left(  1-q\right)  ^{k}}%
\]
is exponentially close to $1$ with respect to $k\rightarrow\infty$ and the
corresponding element in the Green's function $\mathbf{M}_{\mathrm{Green}}$
tends exponentially towards infinity (cf. \eqref{choiceofsq},
\eqref{Estoddrows}). For $n=2k$, $m=k$ and our choices of $q_{\ell}$ and
$\sigma_{\ell}$, it holds%
\begin{equation}
\left\vert \left(  \mathbf{M}_{\mathrm{Green}}^{\left(  2n\right)  }\right)
_{2k+1,2n}\right\vert =\frac{1}{\sqrt{1-|Q_{k}|^{2}}}=\frac{1}{2}\left(
\left(  \frac{1+q}{1-q}\right)  ^{k/2}+\left(  \frac{1-q}{1+q}\right)
^{k/2}\right)  ,\label{equalMgreen}%
\end{equation}
i.e., this matrix entry grows exponentially with respect to increasing $k$. By
adjusting the mesh widths $h_{j}$ and the constant $c$ such that all
$\exp\left(  -2\frac{\operatorname*{i}\omega h_{j+1}}{c_{j+1}}\right)  $
coincide with the phase factors $\sigma_{j}$ as defined in this example, we
have a configuration, where the stability constant grows exponentially in the
number of jumps. In order to achieve
\[
\sigma_{j}=\left\{
\begin{array}
[c]{ll}%
1 & j=k,\\
-1 & \text{otherwise,}%
\end{array}
\right.
\]
we obtain the following conditions for the mesh sizes
\begin{equation}
h_{j}=\frac{\pi}{\omega}\times\left\{
\begin{array}
[c]{ll}%
c_{j}m_{j} & j=k+1,\\
& \\
c_{j}\left(  m_{j}+\frac{1}{2}\right)   & \text{otherwise,}%
\end{array}
\right.  \label{choiceofh}%
\end{equation}
for any sequence $m_{j}\in\mathbb{Z}$. Moreover, the side condition%
\begin{equation}
\sum_{j=1}^{n+1}h_{j}=2\label{sidecond}%
\end{equation}
applies. In (\ref{choiceofh}), we choose $m_{j}=1$ for $j=k+1$ and $m_{j}=0$
otherwise. Hence (\ref{sidecond}) is equivalent to%
\[
2c_{k+1}+\sum_{\substack{i=1\\i\neq k+1}}^{2k+1}c_{i}=\frac{4\omega}{\pi}%
\]
Let $k$ be even. The definition (\ref{defwavespeedex}) implies%
\[
c=\frac{2\omega}{\pi\left(  1-q+k\right)  }.
\]

We summarise the findings of this example. Let $\omega>0$ be given and let
$n=2k $ for even $k\in\mathbb{N}$ and $k= O(\omega)$. Choose some $q\in\left[
0,1\right[  $. Define%
\[
c:=\frac{2\omega}{\pi\left(  1-q+k\right)  } = O(1)
\]
and the piecewise constant wave speed by
\[
c_{j}=\left\{
\begin{array}
[c]{ll}%
c\left(  1-q\right)  & j \text{ odd,}\\
c\left(  1+q\right)  & j \text{ even.}%
\end{array}
\right.
\]
Define mesh sizes according to%
\[
h_{j}=\frac{c_{j}}{\omega}\pi\times\left\{
\begin{array}
[c]{ll}%
1 & j=k+1,\\
\frac{1}{2} & \text{otherwise}%
\end{array}
\right.
\]
and mesh points $x_{j}$ recursively by $x_{0}:=-1$ and for $1\leq j\leq n+1$
by $x_{j}=x_{j-1}+h_{j}$. The stability constant for problem
(\ref{Helmpwconst}) grows exponentially with increasing $\omega$.

\begin{remark}
\label{Remlowerbound}From Lemma \ref{lem:lowerbound} we conclude that the norm
$\Vert u^{\left(  k\right)  }\Vert$ may grow exponentially in $\omega$;
indeed, for this example, the choice $g_{1}=0$, $g_{2}=1$ and $\omega=k$ for
$k\geq2$ leads to%
\[
\frac{4}{3\pi}\leq c\leq\frac{2}{\pi},\text{\quad}\frac{4}{3\pi}\left(
1-q\right)  \leq c_{j}\leq\frac{4}{\pi},\text{\quad}\frac{2}{3}\frac{\left(
1-q\right)  }{\omega}\leq h_{j}\leq\frac{4}{\omega},\text{\quad}\frac{\pi}%
{2}\leq\frac{\omega h_{j}}{c_{j}}\leq\pi.
\]
Thus Lemma \ref{lem:lowerbound} implies%
\[
\left\Vert u^{\left(  k\right)  }\right\Vert _{L^{2}\left(  \Omega\right)
}\geq\frac{2\pi\sqrt{1-q}}{3\sqrt{5}\left(  \pi+2\right)  }\left(  \frac{\pi
}{4}\right)  ^{k}\omega^{k-1/2}\max\left\{  \left\vert A_{j}\right\vert
,\left\vert B_{j}\right\vert \right\}  .
\]
From (\ref{equalMgreen}) and the proof of Theorem \ref{Theo:StabEst} we
conclude that%
\[
\max_{j}\max\left\{  \left\vert A_{j}\right\vert ,\left\vert B_{j}\right\vert
\right\}  \geq\left(  \frac{1+q}{1-q}\right)  ^{\omega/2}\frac{1-q}{3\pi
\omega}%
\]
holds. The combination of these estimates leads to%
\begin{align*}
\left(  \int_{\Omega}|u^{\prime}|^{2}+\left(  \frac{\omega}{c}\right)
^{2}|u|^{2}\right)  ^{\frac{1}{2}} &  \geq\frac{\pi\left(  1-q\right)  ^{3/2}%
}{9\sqrt{5}\left(  \pi+2\right)  }\left(  \frac{1}{2}+\frac{\pi}%
{\sqrt{5\left(  1-q\right)  }\left(  \pi+2\right)  }\right)  \omega
^{-1/2}\left(  \frac{1+q}{1-q}\right)  ^{\omega/2}\\
&  \geq C_{q}\alpha_{q}^{\omega}%
\end{align*}
for some $C_{q}>0$ and $\alpha_{q}\in\left(  0,1\right)  $ depending only on
$q\in\left(  0,1\right)  $.
\end{remark}

\section{Recursive Representation of the Inverse Green's Function}

\label{Sec:inverseGreen} In Section \ref{Sec:generalGrowth}, we discuss the
possible growth of $|Q_{j}|$ for fixed $\omega$ and general parameters $c_{j}$
and $h_{j}$. This will directly lead to a stability estimate if the step sizes
are above resonance, i.e. where $\frac{h}{c}\geq O\left(  \frac{1}{\omega
}\right)  $ (cf. Section \ref{hbounded}). In Section \ref{hnotbounded}, we
will restrict to the case where $c$ oscillates perfectly between two values
and we will show, that for any general configuration of $(h_{j})_{j=1}^{n+1}$,
the stability constant cannot exceed the exponential growth with respect to
$\omega$ described in Section \ref{example_exp}.

\subsection{The Influence of the Phase Factors $\sigma_{i}$}

\label{Sec:generalGrowth}

For the estimate of $Q_{j}$, we have first to introduce some quantities and
conventions. The coefficient $Q_{j}$ depends on $\left(  q_{i}\right)
_{i=1}^{j}$ and $\left(  \sigma_{i}\right)  _{i=1}^{j}$. The index $j$ in
$Q_{j}\left(  \mathbf{q},\mbox{\boldmath$ \sigma$}\right)  $ always indicates
the lengths of $\mathbf{q}=\left(  q_{i}\right)  _{i=1}^{j}$ and
$\mbox{\boldmath$ \sigma$}=\left(  \sigma_{i}\right)  _{i=1}^{j}$, where
$q_{i}\in \left[  -q,q\right]$ and $\sigma_{i}\in\mathcal{C}$ for all $1\leq i\leq j$. For
given $\boldsymbol{q}$ and $\boldsymbol{\sigma}$, let
$\widehat{\boldsymbol{\sigma}}=\widehat{\boldsymbol{\sigma}}\left(
\boldsymbol{q},\boldsymbol{\sigma}\right)  =\left(  \hat{\sigma}_{1}%
,\ldots,\hat{\sigma}_{j-1},\hat{\sigma}_{j}\right)  ^{\intercal}$ be defined
by
\begin{equation}%
\begin{split}
\hat{\sigma}_{i} &  :=\operatorname*{sign}\left(  q_{i}q_{i+1}\right)
\in\left\{  -1,1\right\}  \qquad\forall1\leq i\leq j-1,\\
\hat{\sigma}_{j} &  :=\sigma_{j}.
\end{split}
\label{defsigmaihut}%
\end{equation}
In view of Lemma \ref{LemGnm}, we investigate $\left\vert Q_{n-m}\right\vert $
for $m=1,2,\ldots,n$.

\begin{lemma}
\label{estA1Q1} Recall the definitions of $Q_{j}$ and $Q_{j}^{\prime}$ as in \eqref{defamnqmntot}.

\begin{enumerate}
\item[(i)] It holds%
\begin{equation}
\max_{\substack{\sigma_{i}\in\mathcal{C}\\1\leq i\leq j}}\left\vert
Q_{j}\left(  \mathbf{q},\mbox{\boldmath$ \sigma$}\right)  \right\vert
=\left\vert Q_{j}\left(  \mathbf{q},\widehat{\boldsymbol{\sigma}}\right)
\right\vert ,\label{estqmm}%
\end{equation}
i.e., the maximizer $\left(  \hat{\sigma}_{i}\right)  _{i=1}^{j-1}$ coincides
with the maximizer for a larger $j^{\prime}>j$ in the first $j-1$ components.

\item[(ii)] The quantity $Q_{j}^{\prime}\left(  \mathbf{q}\right)  :=\sigma
_{j}Q_{j}\left(  \mathbf{q},\widehat{\mbox{\boldmath$ \sigma$}}\right)  $
satisfies $Q_{j}^{\prime}\left(  \mathbf{q}\right)  \in\mathbb{R}$,
$\left\vert Q_{j}^{\prime}\left(  \mathbf{q}\right)  \right\vert <1$ and%
\begin{equation}
\operatorname*{sign}\left(  Q_{j}^{\prime}\left(  \mathbf{q}\right)  \right)
=\operatorname*{sign}\left(  q_{j}\right)  . \label{signqm}%
\end{equation}
For given $\mathbf{q}\in \left[  -q,q\right]^{j}$ and $\widehat{\boldsymbol{\sigma}}=\left(
\hat{\sigma}_{i}\right)  _{i=1}^{j}$ as in (\ref{defsigmaihut}) it holds%
\begin{equation}
\operatorname*{sign}\left(  Q_{\ell}\left(  \mathbf{q}%
,\widehat{\mbox{\boldmath$ \sigma$}}\right)  \right)  =\operatorname*{sign}%
\left(  q_{\ell+1}\right)  \qquad\forall1\leq\ell\leq j-1. \label{signql}%
\end{equation}

\item[(iii)] For any sequence $\mathbf{\tilde{q}}=\left(  \tilde{q}_{i}\right)
_{i=1}^{j}$ such that $\tilde{q}_{i}\in\{-q,q\}$ and corresponding
$\boldsymbol{\hat{\sigma}}$, the maximum over all $q_{i}\in\left[  -q,q\right]$ is given
by
\begin{equation}%
\begin{split}
\max_{\mathbf{q}\in \left[  -q,q\right]^{j}}\max_{\substack{\sigma_{i}\in\mathcal{C}\\1\leq
i\leq j}}\left\vert Q_{j}\left(  \mathbf{q},\boldsymbol{\sigma}\right)
\right\vert  &  =\left\vert Q_{j}\left(  \mathbf{\tilde{q}}%
,\widehat{\boldsymbol{\sigma}}(\mathbf{\tilde{q}})\right)  \right\vert \\
&  =\frac{\left(  1+q\right)  ^{j}-\left(  1-q\right)  ^{j}}{\left(
1+q\right)  ^{j}+\left(  1-q\right)  ^{j}}.
\end{split}
\label{fullestQm}%
\end{equation}

\end{enumerate}
\end{lemma}
\proof We prove the equalities in \eqref{estqmm} and \eqref{signqm} by
induction. Both equalities are trivial for $j=1$. For $j=2$, we get%
\[
\max_{\sigma_{1},\sigma_{2}\in\mathcal{C}}\left\vert Q_{2}\left(  \mathbf{q}%
,\mbox{\boldmath$ \sigma$}\right)  \right\vert =\max_{\sigma_{1}\in
\mathcal{C}}\left\vert \frac{q_{2}+\frac{q_{1}}{\sigma_{1}}}{1+\text{$q_{2}$%
}\frac{\text{$q_{1}$}}{\sigma_{1}}}\right\vert =\max_{\sigma_{1}\in
\mathcal{C}}\sqrt{\frac{q_{2}^{2}+q_{1}^{2}+2q_{1}q_{2}\operatorname{Re}%
\sigma_{1}}{1+q_{2}^{2}q_{1}^{2}+2q_{1}q_{2}\operatorname{Re}\sigma_{1}}%
}=\frac{\left\vert q_{1}\right\vert +\left\vert q_{2}\right\vert
}{1+\left\vert q_{1}q_{2}\right\vert }%
\]
and the maximum is achieved for $\sigma_{1}:=\hat{\sigma}_{1}$ with
$\hat{\sigma}_{1}:=\operatorname*{sign}\left(  q_{1}q_{2}\right)  $. In this
case
\[
\operatorname{sign}Q_{2}^{\prime}\left(  \mathbf{q}\right)
=\operatorname{sign}\frac{q_{2}+\frac{q_{1}}{\left(  \operatorname{sign}%
q_{1}q_{2}\right)  }}{1+\frac{q_{2}q_{1}}{\left(  \operatorname{sign}%
q_{1}q_{2}\right)  }}=\operatorname{sign}q_{2}.
\]
Next we consider the case $j\geq3$. By induction we assume that (\ref{estqmm})
and (\ref{signqm}) hold for $j^{\prime}=1,2,\ldots,j-1$. Thus, we get by the
maximum modulus principle%
%TCIMACRO{\TeXButton{maxsigma}{\begin{subequations}
%\label{maxsigma}
%\end{subequations}} }%
%BeginExpansion
\begin{subequations}
\label{maxsigma}
\end{subequations}
%EndExpansion
\begin{equation}
\max_{\substack{\sigma_{i}\in\mathcal{C}\\1\leq i\leq j}}\left\vert
Q_{j}\left(  \mathbf{q},\mbox{\boldmath$ \sigma$}\right)  \right\vert
=\max_{\sigma_{j-1}\in\mathcal{C}}\left\vert \frac{q_{j}+\frac{Q_{j-1}%
^{\prime}\left(  \mathbf{q}\right)  }{\sigma_{j-1}}}{1+q_{j}\frac
{Q_{j-1}^{\prime}\left(  \mathbf{q}\right)  }{\sigma_{j-1}}}\right\vert .
\tag{%
%TCIMACRO{\TeXButton{maxsigma}{\ref{maxsigma}}}%
%BeginExpansion
\ref{maxsigma}%
%EndExpansion
a}\label{maxsigmaa}%
\end{equation}
Similar as before we see that the right-hand side attains its maximum for
$\sigma_{j-1}=\hat{\sigma}_{j-1}$ with $\hat{\sigma}_{j-1}=\operatorname{sign}%
\left(  q_{j}Q_{j-1}^{\prime}\left(  \mathbf{q}\right)  \right)  $, i.e.,%
\begin{equation}
\underset{_{\substack{\sigma_{i}\in\mathcal{C}\\1\leq i\leq j}}}{\max
}\left\vert Q_{j}\right\vert =\frac{\left\vert q_{j+1}\right\vert +\left\vert
Q_{j-1}^{\prime}\left(  \mathbf{q}\right)  \right\vert }{1+\left\vert
q_{j+1}\right\vert \left\vert Q_{j-1}^{\prime}\left(  \mathbf{q}\right)
\right\vert }. \tag{
%TCIMACRO{\TeXButton{maxsigma}{\ref{maxsigma}}}%
%BeginExpansion
\ref{maxsigma}%
%EndExpansion
b}\label{maxsigmab}%
\end{equation}
By induction we have $\operatorname{sign}Q_{j-1}^{\prime}\left(
\mathbf{q}\right)  =\operatorname*{sign}q_{j}$ so that $\hat{\sigma}%
_{j}=\operatorname{sign}\left(  q_{j}q_{j+1}\right)  $. Thus,%
\begin{align*}
\operatorname{sign}\left(  Q_{j}^{\prime}\left(  \mathbf{q}%
,\widehat{\mbox{\boldmath$ \sigma$}}\right)  \right)   &  =\operatorname{sign}%
\left(  \frac{q_{j}+\frac{Q_{j-1}^{\prime}\left(  \mathbf{q}\right)
}{\operatorname{sign}\left(  q_{j-1}q_{j}\right)  }}{1+q_{j}\frac
{Q_{j-1}^{\prime}\left(  \mathbf{q}\right)  }{\hat{\sigma}_{j-1}}}\right)
=\operatorname{sign}\left(  q_{j}+\left(  \operatorname{sign}q_{j}\right)
\left\vert Q_{j-1}^{\prime}\left(  \mathbf{q}\right)  \right\vert \right) \\
&  =\left(  \operatorname{sign}q_{j}\right)  \operatorname{sign}\left(
\left\vert q_{j}\right\vert +\left\vert Q_{j-1}^{\prime}\left(  \mathbf{q}%
\right)  \right\vert \right)  =\operatorname{sign}q_{j}\text{.}%
\end{align*}
The sign of $Q_{j}$ as in (\ref{signql}) can be determined by%
\[
\operatorname*{sign}\left(  Q_{j}\left(  \mathbf{q}%
,\widehat{\mbox{\boldmath$ \sigma$}}\right)  \right)  =\operatorname*{sign}%
\left(  \hat{\sigma}_{j}Q_{j}^{\prime}\right)  =\left(  \operatorname*{sign}%
Q_{j}^{\prime}\right)  \operatorname{sign}\left(  q_{j}q_{j+1}\right)
=\operatorname*{sign}q_{j+1}\text{.}%
\]
Hence part (i) and (ii) are proved.

To prove the bound (\ref{fullestQm}) we observe that the coefficients
$\rho_{j}:=|Q_{j}|$ are majorized by the sequence%
\[
\tilde{r}_{1}\left(  q\right)  =q\text{\quad and }\forall j\geq2\quad\tilde
{r}_{j}:=\frac{q+\tilde{r}_{j-1}}{1+q\tilde{r}_{j-1}}.
\]
The closed form of this recursion is given by%
\begin{equation}
\tilde{r}_{j}\left(  q\right)  =\frac{\left(  1+q\right)  ^{j}-\left(
1-q\right)  ^{j}}{\left(  1+q\right)  ^{j}+\left(  1-q\right)  ^{j}}.
\label{bound:exponentialincrease}%
\end{equation}

\endproof

\begin{remark}
The combination of (\ref{estuk}) and Lemma \ref{estA1Q1}(iii) implies%
\begin{align*}
\left(  \int_{\Omega}|u^{\prime}|^{2}+\left(  \frac{\omega}{c}\right)
^{2}|u|^{2}\right)  ^{\frac{1}{2}}  &  \leq8\frac{c_{\max}}{c_{\min}}%
\max\left\{  \left\vert g_{1}\right\vert ,\left\vert g_{2}\right\vert
\right\}  \max_{1\leq j\leq n}\frac{1}{\sqrt{1-\left\vert Q_{j}\right\vert
^{2}}}\\
&  \leq4\frac{c_{\max}}{c_{\min}}\left(  \kappa^{n/2}+\kappa^{-n/2}\right)
\max\left\{  \left\vert g_{1}\right\vert ,\left\vert g_{2}\right\vert
\right\}
\end{align*}
with the \emph{condition number} $\kappa$ of the wave speed as in
(\ref{condnumber_wavespeed}). This shows that for \emph{fixed} number of jumps the
stability estimate is independent of the wave number $\omega$. Such types of
estimates are also proved in \cite{ChaumontFreletDiss},
\cite{GrahamSauter2017}.
\end{remark}

From now on, we will consider $c$ piecewise constant and perfectly oscillating
between two values $c_{\min}$ and $c_{\max}$, i.e.
\[
c_{j}=%
\begin{cases}
c_{0} & \text{ if }j\text{ is odd, }\\
c_{1} & \text{ if }j\text{ is even, }%
\end{cases}
\]
with $0<c_{\min}=\min\{c_{0},c_{1}\}\leq\max\{c_{0},c_{1}\}=c_{\max}<\infty$.
In that case, we know from Lemma \ref{estA1Q1} that $|Q_{j}|$ increases
maximally for a certain specific choice of $(\sigma_{i})_{i}$. However,
motivated by the example in Section \ref{example_osc}, we know that this can
be very pessimistic. In particular if $h\leq\frac{\delta}{\omega}$ for
sufficiently small $\delta$, then $\sigma_{i-1}=\exp\left(
-2\operatorname*{i}\frac{h_{i}\omega}{c_{i}}\right)  \approx1$. If $\sigma
_{i}=1$ for all $i$ (and $q_{i}=(-1)^{i+1}q$), we have seen that $(Q_{i})_{i}$
can be bounded away from 1 independent of the number of jumps (cf. Section
\ref{example_osc}). The idea we follow is to split the domain into two types
of subsequences of $h_{i}$. The first type covers parts of the domain where
$\omega\frac{h_{i}}{c_{i}}$ is bounded from below. In this case, we can be
bound $|Q_{m}|$ from above with respect to the frequency $\omega$ using
estimates of the type \eqref{bound:exponentialincrease}. For parts of the
domain where $\omega\frac{h_{i}}{c_{i}}$ is small, we will use another
approach, but with a similar result. In the end, the two estimates can be
combined by finding an upper bound on $|Q_{j+m}|$ with respect to the value
$|Q_{j}|$ instead of $|Q_{0}|$. The following corollary restates Lemma
\ref{estA1Q1} for this purpose and the proof is a repetition of arguments.

\begin{corollary}
\label{Cor:maxGrowth} Define the sequence
\[
Q_{j}=\frac{q_{j}+Q_{j-1}}{\sigma_{j}(1+q_{j}Q_{j-1})}%
\]
for $Q_{1}=\frac{\tilde{q}}{\sigma_{1}}$ for $0\leq\tilde{q}<1$. Also assume
that $q_{j}=(-1)^{j+1}q$ for $2\leq j\leq n-1$ and some $0<q<1$. We define
\[
r_{\tilde{q},j}(q):=\frac{(1+\tilde{q})(1+q)^{j-1}-(1-\tilde{q})(1-q)^{j-1}%
}{(1+\tilde{q})(1+q)^{j-1}+(1-\tilde{q})(1-q)^{j-1}}.
\]

\begin{enumerate}
\item[(i)] $r_{\tilde{q},j}$ is increasing in $\tilde{q}$.

\item[(ii)] If $\sigma_{j}=-1$ for all $1\leq j\leq n$, then
\[
|Q_{j}|=r_{\tilde{q},j}(q)
\]
and $\operatorname*{sign}(Q_{j})=(-1)^{j+1}$ for all $1\leq j\leq n$.

\item[(iii)] If $(\sigma_{j})_{j=1}^{\infty}\subset\mathcal{C}$ is a general
sequence, then
\[
|Q_{m+j}|\leq r_{|Q_{j}|,m}(q),\quad\forall~1\leq m+j\leq n.
\]

\end{enumerate}
\end{corollary}

\subsection{Estimate of $Q_{m}$ for Step Sizes above Resonance Case}

\label{hbounded} Assume that $\omega\frac{h_{j}}{c_{j}}>\varepsilon$ for some
$\varepsilon>0$ and for all $1\leq j\leq n+1$. Then we have
\[
\left(  n+1\right)  \cdot\varepsilon\leq\sum_{j=1}^{n+1}\omega\frac{h_{j}%
}{c_{j}}\leq2\frac{\omega}{c_{\min}},
\]
and therefore the length of the sequence $(Q_{j})_{j=1}^{n}$ is bounded by
$n\leq\frac{2\omega}{\varepsilon c_{\min}}-1=N({\varepsilon,\omega})$. Using
Lemma \ref{estA1Q1} and Corollary \ref{Cor:maxGrowth}, we derive
\begin{align}
\label{eq:maxgrowth}|Q_{j}|\leq r_{|Q_{0}|, N({\varepsilon,\omega})}(q),
\quad1\leq j\leq n.
\end{align}
The result obtained so far is summarised in the next theorem.

\begin{proposition}
\label{Prop:aboveresonance} Assume that the step size $h_{i}$ is bounded from
below, i.e. the number of jumps is bounded from above by $\alpha\omega$ for
some $\alpha$, and $c$ is piecewise constant. Then the stability constant in
(\ref{conjecture}) of the Helmholtz problem \eqref{Helmpwconst} satisfies
\[
C_{\operatorname{stab}}\lesssim\alpha_{q}^{-\frac{\omega}{c_{\min}}},~\text{
for some }0<\alpha_{q}<1.
\]

\end{proposition}

This result coincides with the stability estimate proved in
\cite{ChaumontFreletDiss} (for slightly different boundary conditions) and
shows that our theory recovers this result.

\subsection{Estimate of $Q_{m}$ for small Step Sizes}

\label{hnotbounded} In this section, we discuss the growth of the magnitude of
$Q_{j}$ for small widths $h_{j}$, when the wave speed $c$ is perfectly
oscillating between two values. In this case, we cannot control the number of
jumps, however we know that the sum of the widths equals the interval length
$|\Omega|=2$, i.e.
\[
\sum_{j=1}^{n+1}\omega\frac{h_{j}}{c_{j}}\leq\frac{\omega}{c_{\min}}\sum
_{j=1}^{n+1}h_{j}\leq\frac{2\omega}{c_{\min}}%
\]
and this will play a key role in the proof of Proposition \ref{Prop:hsmall}.

\begin{proposition}
\label{Prop:hsmall} We define
\[
q_{j}=\left(  -1\right)  ^{j +1}q
\]
and consider $\sigma_{j}=\operatorname*{e}^{\operatorname*{i}\phi_{j+1}}$ for
$\phi_{j }=-2\omega\frac{h_{j}}{c_{j}}\in\left[  -\phi,0\right]  $ for
sufficiently small $\phi>0$. We set%
\[
Q_{j}:=%
\begin{cases}
Q_{1} & j=1,\\
\dfrac{\left(  -1\right)  ^{j+1}q+Q_{j-1}}{\sigma_{j}\left(  1+\left(
-1\right)  ^{j+1}qQ_{j-1}\right)  } & 2\leq j\leq n,
\end{cases}
\]
for some $Q_{1}$ with $|Q_{1}| <1$. Then the estimate
\begin{align}
\label{boundfinal:hsmall}|Q_{j}|^{2} \leq1-(1-|Q_{1}|^{2})\alpha_{q}%
^{\omega/c_{\min}}%
\end{align}
holds for some $\alpha_{q}\in(0,1)$.
\end{proposition}

\begin{remark}
$\alpha_{q} \rightarrow0$ as $q\rightarrow1$.
\end{remark}

\proof We recall $Q_{j}^{\prime}=\sigma_{j}Q_{j}$, $\rho_{j}:=\left\vert
Q_{j}^{\prime}\right\vert =\left\vert Q_{j}\right\vert $ and $Q_{j}^{\prime
}=:s_{j}\rho_{j}$ for some $s_{j}=\operatorname*{e}^{\operatorname*{i}\psi
_{j}}$ with $\psi_{j}\in\left[  -\pi,\pi\right[  $. We double the iteration and
arrive at
\[
s_{j+2}\rho_{j+2}=\dfrac{\left(  1-\sigma_{j+1}\right)  \left(  -1\right)
^{j+1}q+\left(  1-\sigma_{j+1}q^{2}\right)  \frac{s_{j}}{\sigma_{j}}\rho_{j}%
}{\sigma_{j+1}-q^{2}+\left(  1-\sigma_{j+1}\right)  \left(  -1\right)
^{j}q\frac{s_{j}}{\sigma_{j}}\rho_{j}}.
\]
Two dimensional Taylor expansion w.r.t $(\phi_{j},\phi_{j+1})$ at 0 yields
\begin{equation}
\rho_{j+2}^{2}=\rho_{j}^{2}+\left(  -1\right)  ^{j}\frac{2q}{1-q^{2}}%
(1-\rho_{j}^{2})\rho_{j}\left(  \sin\psi_{j}\right)  \phi_{j+1}+\text{h.o.t.}
\label{recursrho}%
\end{equation}
First, we leave away the higher other term (\textquotedblleft
h.o.t.\textquotedblright) in \eqref{recursrho}, restrict to the case of even
$j$ in \eqref{recursrho}, and study the recursion
\[
\tilde{\rho}_{j+1}^{2}=\tilde{\rho}_{j}^{2}+\frac{2q}{1-q^{2}}\left(
1-\tilde{\rho}_{j}^{2}\right)  \underbrace{\left(  \sin\psi_{2j}\right)
\phi_{2j+1}\tilde{\rho}_{j}}_{=:\gamma_{j}},
\]
for $j\geq1$. This recursion can be resolved and we obtain%
\begin{equation}
\tilde{\rho}_{j+1}^{2}=\tilde{\rho}_{j+1}^{2}\left(  \left(  \gamma
_{k}\right)  _{k=1}^{j}\right)  =\tilde{\rho}_{1}^{2}+\frac{2q}{1-q^{2}%
}\left(  1-\tilde{\rho}_{1}^{2}\right)  \sum_{\ell=1}^{j}\gamma_{\ell}%
\prod\limits_{k=1}^{\ell-1}\left(  1-\frac{2q}{1-q^{2}}\gamma_{k}\right)  ,
\label{fullres}%
\end{equation}
where $\gamma_{k}\in\lbrack-\phi,\phi]$, $\frac{2q}{1-q^{2}}\in\mathbb{R}_{+}$
and $\phi>0$. Now for $\phi\leq\frac{1-q^{2}}{4q}$ and $i\leq j$, we have
\begin{align*}
\frac{\partial\tilde{\rho}_{j+1}^{2}}{\partial\gamma_{i}}=  &  \frac
{2q}{1-q^{2}}(1-\rho_{1}^{2})\frac{\partial}{\partial\gamma_{i}}\left(
\gamma_{i}\prod_{k=1}^{i-1}\left(  1-\frac{2q}{1-q^{2}}\gamma_{k}\right)
+\sum_{\ell=i+1}^{j}\gamma_{\ell}\prod_{k=1}^{\ell-1}\left(  1-\frac
{2q}{1-q^{2}}\gamma_{k}\right)  \right) \\
=  &  \frac{2q}{1-q^{2}}(1-\rho_{1}^{2})\prod_{k=1}^{i-1}\left(  1-\frac
{2q}{1-q^{2}}\gamma_{k}\right)  \left(  1-\frac{2q}{1-q^{2}}\sum_{\ell
=i+1}^{j}\gamma_{\ell}\prod_{k=i+1}^{\ell-1}\left(  1-\frac{2q}{1-q^{2}}%
\gamma_{k}\right)  \right) \\
=  &  \frac{2q}{1-q^{2}}(1-\rho_{1}^{2})\prod_{\substack{k=1\\k\neq i}%
}^{j}\left(  1-\frac{2q}{1-q^{2}}\gamma_{k}\right)  >0.
\end{align*}
This means that the r.h.s of \eqref{fullres} is increasing in $\gamma_{i}$. We
want to find the maximum of \eqref{fullres} under the restriction $\sum
_{i=1}^{j}|\gamma_{i}|<s$, for some $s\in\mathbb{R}_{+}$. If $s\geq n\phi$,
then clearly $\gamma_{i}=\phi$ for all $i$ is the maximizer (see Figure
\ref{Fig:s>nphi}). Since \eqref{fullres} is increasing in every variable, we
can assume $\gamma_{i}\geq0$ and consider the restriction $\sum_{i=1}%
^{j}\gamma_{i}=s$. Let $s<n\phi$, We consider the $(j-1)$-dimensional plane
going through the points $se_{1},se_{2},\dotsc se_{j}$, parametrized by
$\tau:\mathbb{R}^{j-1}\rightarrow\mathbb{R}^{j}$, $(x_{1},\dotsc
,x_{j-1})\rightarrow s(e_{j}+(e_{1}-e_{j})x_{1}+(e_{2}-e_{j})x_{2}%
\dots+(e_{j-1}-e_{j})x_{j-1})$ (see Figure \ref{Fig:L1Linftyinside}). Then
$\tilde{\rho}_{j}^{2}\circ\tau$ is a function defined on $\mathbb{R}^{j-1}$
by
\begin{align*}
\tilde{\rho}_{j+1}^{2}\circ\tau~(\tilde{\gamma}_{1},\tilde{\gamma}_{2}%
,\dotsc,\tilde{\gamma}_{j-1})  &  =\tilde{\rho}_{j+1}^{2}\left(  s\tilde{\gamma
}_{1},\dotsc,s\tilde{\gamma}_{j-1},s\left(  1-\sum_{i=1}^{j-1}\tilde{\gamma
}_{i}\right)  \right) \\
&  =\tilde{\rho}_{1}^{2}+\frac{2q}{1-q^{2}}(1-\tilde{\rho}_{1}^{2})s\sum
_{l=1}^{j-1}\tilde{\gamma}_{l}\prod_{k=1}^{l-1}\left(  1-\frac{2q}{1-q^{2}}%
s\tilde{\gamma}_{k}\right) \\
&  +\frac{2q}{1-q^{2}}(1-\tilde{\rho}_{1}^{2})s\left(  1-\sum_{i=1}%
^{j-1}\tilde{\gamma}_{i}\right)  \prod_{k=1}^{j-1}\left(  1-\frac{2q}{1-q^{2}%
}s\tilde{\gamma}_{k}\right).
\end{align*}
The partial derivatives of $\tilde{\rho}_{j+1}\circ\tau$ are given by
\[
\frac{\partial}{\partial\tilde{\gamma}_{i}}\tilde{\rho}_{j+1}^{2}\circ
\tau=\left(  \frac{2q}{1-q^{2}}\right)  ^{2}s^2(1-\tilde{\rho}_{1}^{2}%
)\prod_{\substack{k=1\\k\neq i}}^{j-1}\left(  1-\frac{2q}{1-q^{2}}s%
\tilde{\gamma}_{k}\right)  \left(  \tilde{\gamma}_{i}-\left(  1-\sum
_{k=1}^{j-1}\tilde{\gamma}_{k}\right)  \right)  .
\]
The r.h.s. is zero for all $i$ iff $\tilde{\gamma}_{k}=\frac{1}{j}$ for all
$k$. Therefore for $(\gamma_{k})_{k=1}^{j}=\tau((\tilde{\gamma}_{k}%
)_{k=1}^{j-1})$, we have $\gamma_{k}=\frac{s}{j}$ for all $k=1,2,\dots,j$, and
$\frac{s}{j}\leq\frac{1-q^{2}}{4q}$, since we also have $\gamma_{k}\leq
\phi\leq\frac{1-q^{2}}{4q}$. In this case, we compute
\[%
\begin{split}
\tilde{\rho}_{j+1}^{2}  &  =\rho_{1}^{2}+\frac{2q}{1-q^{2}}(1-\tilde{\rho}%
_{1}^{2})\sum_{\ell=1}^{j}\frac{s}{j}\left(  1-\frac{2q}{1-q^{2}}\frac{s}%
{j}\right)  ^{l-1}\\
&  =\tilde{\rho}_{1}^{2}+\frac{2q}{1-q^{2}}(1-\tilde{\rho}_{1}^{2})\frac{s}%
{j}\left(  \frac{1-\left(  1-\frac{2q}{1-q^{2}}\frac{s}{j}\right)  ^{j}}%
{\frac{2q}{1-q^{2}}\frac{1}{j}}\right) \\
&  =\tilde{\rho}_{1}^{2}+(1-\tilde{\rho}_{1}^{2})\left(  1-\left(  1-\frac
{2q}{1-q^{2}}\frac{s}{j}\right)  ^{j}\right) \\
&  \leq\tilde{\rho}_{1}^{2}+(1-\tilde{\rho}_{1}^{2})\left(  1-\left(  \frac
{1}{2}\right)  ^{\frac{4qs}{1-q^{2}}}\right) \\
&  \overset{s\leq\frac{4\omega}{c_{\min}}}{\leq}\tilde{\rho}_{1}^{2}%
+(1-\tilde{\rho}_{1}^{2})\left(  1-\alpha_{q}^{\frac{\omega}{c_{\min}}}\right)
\\
&  =1-(1-\tilde{\rho}_{1}^{2})\alpha_{q}^{\frac{\omega}{c_{\min}}},
\end{split}
\]
for some $0<\alpha_{q}<1$. This gives an upper bound on $\rho_{j}$ as stated
in \eqref{boundfinal:hsmall}. However, the maximizer of $\tilde{\rho}_{j}%
\circ\tau$ might also lie on the boundary of the set we are considering. We
investigate this in more details. Consider $\tilde{\rho}_{j}\circ\tau$
restricted to $D\subset\mathbb{R}^{j-1}$ such that $\tau(D)$ is the the
intersection of the cube $(0,\phi]^{j}$ and the set $\{x\in\mathbb{R}^{j}%
|\sum_{i=1}^{j}x_{i}=s\}$ (see Figure \ref{Fig:L1Linftyboundary}).
\begin{figure}[ptb]
\begin{subfigure}[t]{0.31\textwidth}
\begin{center}
\begin{tikzpicture}[line cap=round,line join=round,x=0.5cm,y=0.5cm,scale=2]
\draw[<->,color=black] (-1.5,0) -- (1.5,0);
\draw[<->,color=black] (0,-1.5) -- (0,1.5);
\draw(-0.5,0.5) -- (0.5,0.5)-- (0.5,-0.5)--(-0.5,-0.5)--cycle;
\draw (0,1.3) -- (1.3,0) -- (0,-1.3) -- (-1.3,0) -- cycle;
\fill[color = red] (0.5,0.5)circle[radius=1pt];
\end{tikzpicture}
\end{center}
\caption{If \(s\geq n\phi\), the maximum is achieved at the vertex \tikz\draw[red,fill=red] (0,0) circle (.4ex); . }\label{Fig:s>nphi}
\end{subfigure}
\hfill\begin{subfigure}[t]{0.31\textwidth}
\begin{center}
\begin{tikzpicture}[line cap=round,line join=round,x=0.5cm,y=0.5cm,scale=2]
\draw[<->,color=black] (-1.5,0) -- (1.5,0);
\draw[<->,color=black] (0,-1.5) -- (0,1.5);
\draw(-1,1) -- (1,1)-- (1,-1)--(-1,-1)--cycle;
\draw (0,1.3) -- (1.3,0) -- (0,-1.3) -- (-1.3,0) -- cycle;
%\draw [dashed, color = red, line width=1pt] (-0.5,1.8)--(1.8,-0.5);
\draw [dashed, color = red, line width=2pt] (1,0.3) -- (0.3,1);
\end{tikzpicture}
\end{center}
\caption{In the first step, we look for a maximizer along the dashed line.}\label{Fig:L1Linftyinside}
\end{subfigure}
\hfill\begin{subfigure}[t]{0.31\textwidth}
\begin{center}
\begin{tikzpicture}[line cap=round,line join=round,x=0.5cm,y=0.5cm,scale=2]
\draw[<->,color=black] (-1.5,0) -- (1.5,0);
\draw[<->,color=black] (0,-1.5) -- (0,1.5);
\draw(-1,1) -- (1,1)-- (1,-1)--(-1,-1)--cycle;
\draw (0,1.3) -- (1.3,0) -- (0,-1.3) -- (-1.3,0) -- cycle;
\fill[color = red] (0.3,1)circle[radius=1pt];
\fill[color = red] (1,0.3)circle[radius=1pt];
\end{tikzpicture}
\end{center}
\caption{In the second step, we check if the points at the boundary (\tikz\draw[red,fill=red] (0,0) circle (.4ex);) are maximizers.}\label{Fig:L1Linftyboundary}
\end{subfigure}
\caption{Showing the different situations for finding the maximum of
\eqref{fullres} under the restriction $\sum_{i=1}^{j}|\gamma_{i}|\leq s$.}%
\end{figure}Now we note that on the boundary of this set we have $\gamma
_{k}\in\{0,\phi\}$ for some $k$. Finding the maximum on the boundary can be
reduced to the same problem of dimension $n-2$ (fixing $\gamma_{k}\in
\{0,\phi\}$) with $\tilde{s}=s-\gamma_{k}$. Since $\tilde{\rho}_{j}$ is
symmetric in all its variables we find that for $m$ s.t. $(m-1)\phi<s\leq
m\phi$ a maximizer is given by
\begin{align*}
\gamma_{k}  &  =\phi,\qquad\quad\text{for }k=1,\ldots,m-1\\
\gamma_{m}  &  =s-(m-1)\phi,\\
\gamma_{k}  &  =0,\qquad\quad\text{for }k>m.
\end{align*}
With this choice of $\gamma_{k}$, we have
\begin{align}\label{eq:estpj}
\begin{split}
\tilde{\rho}_{j+1}^{2}  &  \leq\tilde{\rho}_{1}^{2}+\frac{2q}{1-q^{2}}\left(
1-\tilde{\rho}_{1}^{2}\right)  \sum_{l=1}^{m}\gamma_{l}\prod_{k=1}%
^{l-1}\left(  1-\frac{2q}{1-q^{2}}\gamma_{k}\right) \\
& \leq\tilde{\rho}_{1}%
^{2}+\frac{2q}{1-q^{2}}(1-\tilde{\rho}_{1}^{2})\phi\sum_{l=1}^{m}\left(
1-\frac{2q}{1-q^{2}}\phi\right)  ^{l-1}\\
&  =\tilde{\rho}_{1}^{2}+\frac{2q}{1-q^{2}}\left(  1-\tilde{\rho}_{1}%
^{2}\right)  \phi\frac{1-(1-\frac{2q}{1-q^{2}}\phi)^{m}}{\frac{2q}{1-q^{2}%
}\phi}=\tilde{\rho}_{1}^{2}+\left(  1-\tilde{\rho}_{1}^{2}\right)  \left(
1-\left(  1-\frac{2q}{1-q^{2}}\phi\right)  ^{m}\right) \\
&  \leq\tilde{\rho}_{1}^{2}+(1-\tilde{\rho}_{1}^{2})\left(  1-\left(
1-\frac{2q}{1-q^{2}}\phi\right)  ^{s/\phi+1}\right)  =1-(1-\tilde{\rho}%
_{1}^{2})\left(  1-\frac{2q}{1-q^{2}}\phi\right)  ^{s/\phi+1}\\
&  \leq1-(1-\tilde{\rho}_{1}^{2})\alpha_{q}^{\omega/c_{\min}},
\end{split}
\end{align}
using $\phi<\frac{1-q}{4q}$ and $s\leq4\omega/c_{\min}$, with a possibly
adjusted $\alpha_{q}\in(0,1)$ from (\ref{eq:estpj}).\newline Now we consider
the general case including the higher other terms in \eqref{recursrho}. Let
$j$ be even. Taylor expansion w.r.t. $(\phi_{j},\phi_{j+1})$ yields
\begin{align*}
\rho_{j+2}^{2}  &  \leq\rho_{j}^{2}+\frac{2q}{1-q^{2}}\left(  1-\rho_{j}%
^{2}\right)  \rho_{j}\left(  \sin\psi_{j}\right)  \phi_{j+1}\\
&  \phantom{=\rho_{j}^2\cdot}+\frac{2q}{1-q^{2}}(1-\rho_{j}^{2})K(\phi
_{j}+\phi_{j+1})^{2}\\
&  \leq\rho_{j}^{2}+\frac{2q}{1-q^{2}}(1-\rho_{j}^{2})(|\phi_{j+1}%
|+K(|\phi_{j}|+|\phi_{j+1}|)^{2}),
\end{align*}
for a constant $K\in\mathbb{R}_{>0}$ uniformly bounded if $\phi<\frac{1}{8}$.
Now we compute
\begin{align*}
\rho_{j+2}^{2}  &  \leq\rho_{j}^{2}+\frac{2q}{1-q^{2}}(1-\rho_{j}^{2}%
)(|\phi_{j+1}|+K(|\phi_{j}|+|\phi_{j+1}|)^{2})\\
&  =\rho_{j}^{2}\left(  1-\frac{2q}{1-q^{2}}(\underbrace{|\phi_{j+1}%
|+K(|\phi_{j}|+|\phi_{j+1}|)^{2}}_{=:\eta_{j}})\right)  +\frac{2q}{1-q^{2}%
}(|\phi_{j+1}|+K(|\phi_{j}|+|\phi_{j+1}|)^{2})\\
&  =\rho_{j}^{2}\left(  1-\frac{2q}{1-q^{2}}\eta_{j}\right)  +\frac
{2q}{1-q^{2}}\eta_{j}.
\end{align*}
For $\phi<\min\left\{  \frac{1}{4K},\frac{1-q^{2}}{4q}\right\}  $, we have
that $(1-\frac{2q}{1-q^{2}}\eta_{j})$ is positive and hence, we can consider
the majorant
\[
\rho_{j+2}^{\prime2}=\rho_{j}^{\prime2}\left(  1-\frac{2q}{1-q^{2}}\eta
_{j}\right)  +\frac{2q}{1-q^{2}}\eta_{j}%
\]
for $0\leq\eta_{j}<\frac{1-q^{2}}{2q}$ and resolve the representation
\[
\rho_{j+1}^{\prime2}=\rho_{1}^{\prime2}+\frac{2q}{1-q^{2}}(1-\rho_{1}%
^{\prime2})\sum_{\ell=1}^{j/2}\gamma_{\ell}\prod_{k=1}^{\ell-1}\left(
1-\frac{2q}{1-q^{2}}\gamma_{k}\right)  ,
\]
which is increasing in $\gamma_{\ell}:=\eta_{2\ell}$. By the same argument as
before we receive
\[
\rho_{j}^{2}\leq1-(1-\rho_{1}^{2})\alpha_{q}^{\omega/c_{\min}},
\]
for some $0<\alpha_{q}<1$. \endproof

\begin{proposition}
\label{Prop:smallstep} Consider the Helmholtz problem \eqref{Helmpwconst} with
diffusion coefficient $a=1$, $f=0$ and boundary values $g_{1}, g_{2}$. Let $c$
be oscillating between $c_{\min}$ and $c_{\max}$. If $h_{j} \lesssim\frac
{1}{\omega}$ is small enough for all $j$ then the stability constant in
\eqref{conjecture} can be bounded independently of the number of jumps by
\begin{align*}
C_{\operatorname{stab}}\leq C_{q} \alpha_{q}^{-\frac{\omega}{c_{\min}}},
\end{align*}
for some $0<\alpha_{q}<1$ and $C_{q}\in\mathbb{R} $.
\end{proposition}

\subsection{Stability Estimate\label{SecStabEst}}

Proposition \ref{Prop:aboveresonance} and \ref{Prop:smallstep} can be combined
to find a final stability estimate \eqref{conjecture} for the Helmholtz problem.

\begin{theorem}
\label{Th:LastTh} Let $c$ be piecewise constant and perfectly oscillating
between two values, i.e.
\[
c_{j}=%
\begin{cases}
c_{0} & \text{ if }j\text{ is odd, }\\
c_{1} & \text{ if }j\text{ is even, }%
\end{cases}
\]
with $0<c_{\min}=\min\{c_{0},c_{1}\}\leq\max\{c_{0},c_{1}\}=c_{\max}<\infty$.
Let $u$ be the (weak) solution of the Helmholtz problem
\begin{equation}
-u^{\prime\prime}-\left(  \frac{\omega}{c}\right)  ^{2}u=0\text{ in }%
\Omega=\left(  -1,1\right)  ,
\end{equation}
with boundary conditions
\begin{align*}
-u^{\prime}-\operatorname*{i}\frac{\omega}{c_{1}}u=g_{1} &  \text{ at }x=-1,\\
u^{\prime}-\operatorname*{i}\frac{\omega}{c_{n}}u=g_{2} &  \text{ at }x=1.
\end{align*}
Then $u$ satisfies
\[
\left(  \int_{\Omega}\left\vert u^{\prime}\right\vert ^{2}+\left(
\frac{\omega}{c}\right)  ^{2}|u|^{2}\right)  ^{\frac{1}{2}}\leq
C_{\operatorname{stab}}\max\{|g_{1}|,|g_{2}|\},
\]
with
\begin{equation}
C_{\operatorname{stab}}\leq C_{q}\alpha_{q}^{-\frac{\omega}{c_{\min}}}%
\quad\text{ for some }0<\alpha_{q}<1\text{ and }C_{q}\in\mathbb{R}%
.\label{Th:FinalEstimate}%
\end{equation}
The estimate (\ref{Th:FinalEstimate}) is independent of the number of jumps of
$c$ and does not require any periodicity of the media.
\end{theorem}

\subsection{Remarks on the Energy Estimate}

We make some remarks on the achieved results.

\begin{enumerate}
\item In the previous section, we derived an estimate on the stability
constant $C_{\operatorname{stab}}$ with constant principal part. However, the
case where $a$ is not constant can be reduced to the constant case by
introducing the new function
\[
u\left(  x\right)  =v\circ\eta\left(  x\right)
\]
with $\eta:\left(  -1,1\right)  \rightarrow\left(  -1,1\right)  $%
\[
\eta\left(  x\right)  =-1+\frac{2}{A}\int_{-1}^{x}\frac{1}{a\left(  s\right)
}ds\quad\text{and\quad}A=\int_{-1}^{1}\frac{ds}{a\left(  s\right)  }.
\]
Then, it is an easy exercise to verify that $v$ satisfies a Helmholtz equation
with the Laplacian as its principle part, again with Robin boundary conditions.

\item The result in Theorem \ref{Th:FinalEstimate} is valid for coefficients
$c$ which oscillate perfectly between two values. With the developed technique
described above the case where $c$ is monotone can be handled similarly to
Section \ref{hbounded} (but independent of $n$). For arbitrary piecewise
configuration of $c$, the handling of the sequence $(Q_{j})$ requires more
technicalities and is an open question.

\item If the configuration of $c$ is \textit{fixed}, one can estimate the
growth of $|Q_{j}|$ with \eqref{eq:maxgrowth} where $n$ is the number of jumps
of $c$, in particular independent of $\omega$ (and $\varepsilon$), to show that
the stability constant is bounded from above independently of $\omega$. This
was also shown in the PhD thesis \cite{ChaumontFreletDiss}.
\end{enumerate}

\section{Proof of the Representation Formulas\label{ProofsRep}}

\textbf{Proof of Lemma \ref{LemLGS}.} We consider the ansatz
(\ref{localhomsol}). The transmission condition leads to
\begin{equation}
\left[
\begin{array}
[c]{rrrrrrr}%
\mathbf{C}^{\left(  1\right)  } & \mathbf{-G}^{\left(  1\right)  } &
\mathbf{0} & \ldots &  &  & \mathbf{0}\\
\mathbf{0} & \mathbf{C}^{\left(  2\right)  } & -\mathbf{G}^{\left(  2\right)
} &  &  &  & \vdots\\
\vdots &  & \ddots & \ddots &  &  & \\
&  &  &  &  &  & \\
&  &  &  &  &  & \mathbf{0}\\
\mathbf{0} &  &  &  &  & \mathbf{C}^{\left(  n\right)  } & -\mathbf{G}%
^{\left(  n\right)  }%
\end{array}
\right]  \left(
\begin{array}
[c]{c}%
A_{1}\\
B_{1}\\
\vdots\\
A_{n+1}\\
B_{n+1}%
\end{array}
\right)  =\mathbf{0}, \label{transcond}%
\end{equation}
where $\mathbf{C}^{\left(  i\right)  },\mathbf{G}^{\left(  i\right)  }%
\in\mathbb{R}^{2\times2}$ are given by%
\begin{align*}
\mathbf{C}^{\left(  i\right)  }  &  :=\left[
\begin{array}
[c]{rr}%
c_{11}^{\left(  i\right)  } & c_{12}^{\left(  i\right)  }\\
c_{21}^{\left(  i\right)  } & c_{22}^{\left(  i\right)  }%
\end{array}
\right]  :=\left[
\begin{array}
[c]{rr}%
\alpha_{i,i} & \dfrac{1}{\alpha_{i,i}}\\
\dfrac{\alpha_{i,i}}{c_{i}} & -\dfrac{1}{\alpha_{i,i}c_{i}}%
\end{array}
\right]  ,\\
\mathbf{G}^{\left(  i\right)  }  &  :=\left[
\begin{array}
[c]{rr}%
d_{11}^{\left(  i\right)  } & d_{12}^{\left(  i\right)  }\\
d_{21}^{\left(  i\right)  } & d_{22}^{\left(  i\right)  }%
\end{array}
\right]  :=\left[
\begin{array}
[c]{rr}%
\alpha_{i+1,i} & \dfrac{1}{\alpha_{i+1,i}}\\
\dfrac{\alpha_{i+1,i}}{c_{i+1}} & -\dfrac{1}{\alpha_{i+1,i}c_{i+1}}%
\end{array}
\right]  .
\end{align*}
Inserting the boundary conditions into (\ref{transcond}) results in the
following block-tridiagonal system%
\begin{equation}
\left[
\begin{array}
[c]{ccccc}%
\mathbf{S}^{\left(  1\right)  } & \mathbf{T}^{\left(  1\right)  } & \mathbf{0}
& \dots & \mathbf{0}\\
\mathbf{R}^{\left(  1\right)  } & \mathbf{S}^{\left(  2\right)  } & \ddots &
\ddots & \vdots\\
\mathbf{0} & \ddots & \ddots &  & \mathbf{0}\\
\vdots & \ddots &  &  & \mathbf{T}^{\left(  n-1\right)  }\\
\mathbf{0} & \dots & \mathbf{0} & \mathbf{R}^{\left(  n-1\right)  } &
\mathbf{S}^{\left(  n\right)  }%
\end{array}
\right]  \left(
\begin{array}
[c]{c}%
B_{1}\\
A_{2}\\
B_{2}\\
\vdots\\
A_{n}\\
B_{n}\\
A_{n+1}%
\end{array}
\right)  =\left[
\begin{array}
[c]{rr}%
-c_{11}^{\left(  1\right)  } & 0\\
-c_{21}^{\left(  1\right)  } & \vdots\\
0 & 0\\
\vdots & d_{12}^{\left(  n\right)  }\\
0 & d_{22}^{\left(  n\right)  }%
\end{array}
\right]  \left(
\begin{array}
[c]{c}%
A_{1}\\
B_{n+1}%
\end{array}
\right)  \label{LGS1}%
\end{equation}
with%
\begin{align*}
\mathbf{R}^{\left(  i\right)  }  &  =\left[
\begin{array}
[c]{cc}%
0 & c_{11}^{\left(  i+1\right)  }\\
0 & c_{21}^{\left(  i+1\right)  }%
\end{array}
\right]  =\left[
\begin{array}
[c]{cc}%
0 & \alpha_{i+1,i+1}\\
0 & \frac{\alpha_{i+1,i+1}}{c_{i+1}}%
\end{array}
\right]  ,\\
\mathbf{S}^{\left(  i\right)  }  &  =\left[
\begin{array}
[c]{cc}%
c_{12}^{\left(  i\right)  } & -d_{11}^{\left(  i\right)  }\\
c_{22}^{\left(  i\right)  } & -d_{21}^{\left(  i\right)  }%
\end{array}
\right]  =\left[
\begin{array}
[c]{cc}%
\frac{1}{\alpha_{i,i}} & -\alpha_{i+1,i}\\
-\frac{1}{\alpha_{i,i}c_{i}} & -\frac{\alpha_{i+1,i}}{c_{i+1}}%
\end{array}
\right]  ,\\
\mathbf{T}^{\left(  i\right)  }  &  =\left[
\begin{array}
[c]{cc}%
-d_{12}^{\left(  i\right)  } & 0\\
-d_{22}^{\left(  i\right)  } & 0
\end{array}
\right]  =\left[
\begin{array}
[c]{cc}%
-\frac{1}{\alpha_{i+1,i}} & 0\\
\frac{1}{\alpha_{i+1,i}c_{i+1}} & 0
\end{array}
\right]  .
\end{align*}

Next we transform (\ref{LGS1}) to a tridiagonal system. Let
$\operatorname*{row}\left(  i\right)  $ denote the $i$-th row of (\ref{LGS1}).
First, we replace for $1\leq i\leq n-1$ the $\operatorname*{row}\left(
2i\right)  $ by $\rho_{i}\operatorname*{row}\left(  2i-1\right)  -c_{i}%
\rho_{i}\operatorname*{row}\left(  2i\right)  $. This leads to%
\begin{equation}
\left[
\begin{array}
[c]{ccccc}%
\mathbf{\tilde{S}}^{\left(  1\right)  } & \mathbf{\tilde{T}}^{\left(
1\right)  } & \mathbf{0} & \dots & \mathbf{0}\\
\mathbf{\tilde{R}}^{\left(  1\right)  } & \mathbf{\tilde{S}}^{\left(
2\right)  } & \ddots & \ddots & \vdots\\
\mathbf{0} & \ddots & \ddots &  & \mathbf{0}\\
\vdots & \ddots &  &  & \mathbf{\tilde{T}}^{\left(  n-1\right)  }\\
\mathbf{0} & \dots & \mathbf{0} & \mathbf{\tilde{R}}^{\left(  n-1\right)  } &
\mathbf{\tilde{S}}^{\left(  n\right)  }%
\end{array}
\right]  \left(
\begin{array}
[c]{c}%
B_{1}\\
A_{2}\\
B_{2}\\
\vdots\\
A_{n}\\
B_{n}\\
A_{n+1}%
\end{array}
\right)  =\left[
\begin{array}
[c]{rr}%
-\alpha_{1,1} & 0\\
0 & \vdots\\
\vdots & 0\\
& \frac{1}{\alpha_{n+1,n}}\\
0 & \frac{\rho_{n}}{\alpha_{n+1,n}}\frac{c_{n}+c_{n+1}}{c_{n+1}}%
\end{array}
\right]  \left(
\begin{array}
[c]{c}%
A_{1}\\
B_{n+1}%
\end{array}
\right)  \label{LGS2}%
\end{equation}
with%
\[
\mathbf{\tilde{R}}^{\left(  i\right)  }=\left[
\begin{array}
[c]{cc}%
0 & \alpha_{i+1,i+1}\\
0 & 0
\end{array}
\right]  ,\quad\mathbf{\tilde{S}}^{\left(  i\right)  }=\left[
\begin{array}
[c]{cc}%
\frac{1}{\alpha_{i,i}} & -\alpha_{i+1,i}\\
2\frac{\rho_{i}}{\alpha_{i,i}} & \rho_{i}\alpha_{i+1,i}\left(  \frac
{c_{i}-c_{i+1}}{c_{i+1}}\right)
\end{array}
\right]  ,\quad\mathbf{\tilde{T}}^{\left(  i\right)  }=\left[
\begin{array}
[c]{cc}%
-\frac{1}{\alpha_{i+1,i}} & 0\\
-\frac{\rho_{i}}{\alpha_{i+1,i}}\frac{c_{i}+c_{i+1}}{c_{i+1}} & 0
\end{array}
\right]  .
\]
In the next step, we replace for $1\leq i\leq n$ the $\operatorname*{row}%
\left(  2i-1\right)  $ by $\rho_{i}\frac{c_{i}+c_{i+1}}{c_{i+1}}\delta
_{i}\operatorname*{row}\left(  2i-1\right)  -\delta_{i}\operatorname*{row}%
\left(  2i\right)  $. This leads to%
\begin{equation}%
\begin{bmatrix}
\mathbf{\hat{S}}^{(1)} & \mathbf{\hat{T}}^{(1)} & \mathbf{0} & \dots &
\mathbf{0}\\
\mathbf{\hat{R}}^{(1)} & \mathbf{\hat{S}}^{(2)} & \ddots & \ddots & \vdots\\
\mathbf{0} & \ddots & \ddots &  & \mathbf{0}\\
\vdots & \ddots &  &  & \mathbf{\hat{T}}^{(n-1)}\\
\mathbf{0} & \dots & \mathbf{0} & \mathbf{\hat{R}}^{(n-1)} & \mathbf{\hat{S}%
}^{(n)}\\
&  &  &  &
\end{bmatrix}
\left(
\begin{array}
[c]{c}%
B_{1}\\
A_{2}\\
B_{2}\\
\vdots\\
A_{n}\\
B_{n}\\
A_{n+1}%
\end{array}
\right)  =\left[
\begin{array}
[c]{cc}%
-\rho_{1}\frac{c_{1}+c_{2}}{c_{2}}\delta_{1}\alpha_{1,1} & 0\\
0 & \vdots\\
\vdots & 0\\
& 0\\
0 & \frac{\rho_{n}}{\alpha_{n+1,n}}\frac{c_{n}+c_{n+1}}{c_{n+1}}%
\end{array}
\right]  \left(
\begin{array}
[c]{c}%
A_{1}\\
B_{n+1}%
\end{array}
\right)  \label{LGS3}%
\end{equation}
with%
\begin{align*}
\mathbf{\hat{R}}^{\left(  i\right)  }  &  =\left[
\begin{array}
[c]{cc}%
0 & \rho_{i+1}\delta_{i+1}\frac{c_{i+1}+c_{i+2}}{c_{i+2}}\alpha_{i+1,i+1}\\
0 & 0
\end{array}
\right]  ,\\
\mathbf{\hat{S}}^{\left(  i\right)  }  &  =\left[
\begin{array}
[c]{cc}%
\dfrac{\rho_{i}\delta_{i}}{\alpha_{i,i}}\dfrac{c_{i}-c_{i+1}}{c_{i+1}} &
-2\dfrac{c_{i}}{c_{i+1}}\delta_{i}\rho_{i}\alpha_{i+1,i}\\
2\frac{\rho_{i}}{\alpha_{i,i}} & \rho_{i}\alpha_{i+1,i}\dfrac{c_{i}-c_{i+1}%
}{c_{i+1}}%
\end{array}
\right]  ,\\
\mathbf{\hat{T}}^{\left(  i\right)  }  &  =\left[
\begin{array}
[c]{cc}%
0 & 0\\
-\dfrac{\rho_{i}}{\alpha_{i+1,i}}\dfrac{c_{i}+c_{i+1}}{c_{i+1}} & 0
\end{array}
\right]  .
\end{align*}
We choose
\[
\delta_{i}=-\frac{1}{\alpha_{i,i}\alpha_{i+1,i}}\frac{c_{i+1}}{c_{i}}%
\quad\text{and\quad}\rho_{i}=\frac{\alpha_{i+1,i}}{c_{i}+c_{i+1}}%
\]
to obtain the following symmetric tridiagonal system of linear equations%
\begin{equation}
\underset{=:\mathbf{\check{M}}^{\left(  2n\right)  }}{\underbrace{\left[
\begin{array}
[c]{ccccc}%
\mathbf{\check{S}}^{\left(  1\right)  } & \mathbf{\check{T}}^{\left(
1\right)  } & \mathbf{0} & \dots & \mathbf{0}\\
\left(  \mathbf{\check{T}}^{\left(  1\right)  }\right)  ^{\intercal} &
\mathbf{\check{S}}^{\left(  2\right)  } & \ddots & \ddots & \vdots\\
\mathbf{0} & \ddots & \ddots &  & \mathbf{0}\\
\vdots & \ddots &  &  & \mathbf{\check{T}}^{\left(  n-1\right)  }\\
\mathbf{0} & \dots & \mathbf{0} & \left(  \mathbf{\check{T}}^{\left(
n-1\right)  }\right)  ^{\intercal} & \mathbf{\check{S}}^{\left(  n\right)  }%
\end{array}
\right]  }}\underset{\mathbf{x}^{\left(  2n\right)  }}{\underbrace{\left(
\begin{array}
[c]{c}%
B_{1}\\
A_{2}\\
B_{2}\\
\vdots\\
A_{n}\\
B_{n}\\
A_{n+1}%
\end{array}
\right)  }}=\underset{\mathbf{r}^{\left(  2n\right)  }}{\underbrace{\left[
\begin{array}
[c]{rr}%
\frac{1}{c_{1}} & 0\\
0 & \vdots\\
\vdots & 0\\
& 0\\
0 & \frac{1}{c_{n+1}}%
\end{array}
\right]  \left(
\begin{array}
[c]{c}%
A_{1}\\
B_{n+1}%
\end{array}
\right)  }} \label{LGS4}%
\end{equation}
with the $2\times2$ matrices
\[
\mathbf{\check{S}}^{\left(  i\right)  }:=\left[
\begin{array}
[c]{cc}%
\frac{c_{i+1}-c_{i}}{c_{i+1}+c_{i}}\frac{1}{c_{i}\alpha_{i,i}^{2}} & \frac
{2}{c_{i}+c_{i+1}}\frac{\alpha_{i+1,i}}{\alpha_{i,i}}\\
\frac{2}{c_{i}+c_{i+1}}\frac{\alpha_{i+1,i}}{\alpha_{i,i}} & -\frac
{c_{i+1}-c_{i}}{c_{i+1}+c_{i}}\frac{\alpha_{i+1,i}^{2}}{c_{i+1}}%
\end{array}
\right]  ,\quad\mathbf{\check{T}}^{\left(  i\right)  }:=\left[
\begin{array}
[c]{cc}%
0 & 0\\
-\frac{1}{c_{i+1}} & 0
\end{array}
\right]  .
\]

It is easy to see that $\left(  \mathbf{\check{M}}^{\left(  2n\right)
}\right)  ^{-1}$ can be factorized $\left(  \mathbf{\check{M}}^{\left(
2n\right)  }\right)  ^{-1}=\mathbf{D}^{\left(  2n\right)  }\mathbf{M}%
_{\mathrm{Green}}^{\left(  2n\right)  }\mathbf{D}^{\left(  2n\right)  }$ with
$\mathbf{D}^{\left(  2n\right)  }$ and $\mathbf{M}_{\mathrm{Green}}^{\left(
2n\right)  }$ as in (\ref{defD}) and (\ref{defM2n-2}).

\endproof
\medskip\newline\textbf{Proof of Lemma \ref{LemDet}.} First, we will prove%
\begin{equation}
\det\mathbf{M}^{\left(  2n\right)  }=\left(  -1\right)  ^{n}\tilde{p}_{n}
\label{rekpn-1}%
\end{equation}
for $n\geq1$ by induction with $\tilde{p}_{n}$ as in (\ref{defpntilde}). It is
easy to check that $\det\mathbf{M}^{\left(  2\right)  }=-1=-\tilde{p}_{1}$.
From now on we assume that (\ref{rekpn-1}) holds for $n^{\prime}=2,3,\ldots
n-1$. From (\ref{detrek}) we derive%
\begin{align*}
\det\mathbf{M}^{\left(  2n\right)  }  &  =-q_{n} \det\mathbf{M}^{\left(
2n-1\right)  }-\left(  1-q_{n} ^{2}\right)  \det\mathbf{M}^{\left(
2n-2\right)  }\\
\det\mathbf{M}^{\left(  2n-1\right)  }  &  =q_{n} \det\mathbf{M}^{\left(
2n-2\right)  }-\frac{1}{\sigma_{n-1} }\det\mathbf{M}^{\left(  2n-3\right)  }.
\end{align*}
and, in turn, we get
\begin{subequations}
\begin{align}
\det\mathbf{M}^{\left(  2n \right)  }  &  =-\det\mathbf{M}^{\left(
2n-2\right)  }+\frac{q_{n }}{\sigma_{n-1}}\sum_{\ell=1}^{n-1}\frac{\left(
-1\right)  ^{\ell+1}q_{n-\ell}}{\prod\nolimits_{k=2}^{\ell}\sigma_{n-k}}%
\det\mathbf{M}^{ \left(  2\left(  n-\ell-1\right)  \right)  }%
,\label{doubleinda}\\
\det\mathbf{M}^{\left(  2n-1\right)  }  &  =\sum_{\ell=1}^{n}\left(
-1\right)  ^{\ell+1}\frac{q_{n+1-\ell}}{ {\textstyle\prod\nolimits_{k=n+1-
\ell}^{n-1}} \sigma_{k}}\det\mathbf{M}^{\left(  2n-2\ell\right)  }.
\label{doubleindb}%
\end{align}

We insert the induction assumption (\ref{rekpn-1}) for $\tilde{p}_{n}$ into
the right-hand side of (\ref{doubleinda}). Then (\ref{rekpn-1}) follows if we
prove%
\end{subequations}
\begin{equation}
\tilde{p}_{n }=\tilde{p}_{n-1}+\frac{q_{n}}{\sigma_{n-1}}\sum_{\ell=1}%
^{n-1}\frac{q_{n-\ell}}{\prod\nolimits_{k=2}^{\ell}\sigma_{n-k}}\tilde
{p}_{n-\ell-1} \label{ptildenrek}%
\end{equation}
for $n\geq1$. It is simple to verify this equality for $n=1$ and we assume for
the following that this holds for $n^{\prime}=1,2,\ldots,n-1$.

For $1\leq k\leq n-1$, we set
\[
\delta_{n,k}:=\sum_{\ell=k}^{n-1}\frac{q_{n-\ell}}{\prod\nolimits_{j=2}^{\ell
}\sigma_{n-j}}\tilde{p}_{n-\ell-1}%
\]
and prove%
\begin{equation}
\delta_{n,k}=\frac{{\textstyle\prod\nolimits_{j=1}^{n-k}} \sigma_{j}}
{{\textstyle\prod\nolimits_{j=1}^{n-2}} \sigma_{j}} \left(
{\displaystyle\prod\limits_{j=1}^{n-1-k}} \left(  1+q_{j+1}Q_{j}\right)
\right)  Q_{n-k} \label{defdeltan}%
\end{equation}
by induction over $k=n-1,n-2,\ldots1$. We denote the right-hand side in
(\ref{defdeltan}) by $\tilde{\delta}_{n,k}$ and prove $\delta_{n,k}%
=\tilde{\delta}_{n,k}$. For $k=n-1$ we have
\[
\tilde{\delta}_{n,n-1}=\frac{q_{1}}{{\textstyle\prod\nolimits_{j=1}^{n-2}}
\sigma_{j}} \quad\delta_{n,n-1}=\frac{q_{1}}{\prod\nolimits_{k=2}^{n-1}%
\sigma_{n-k}}.
\]
Assume the assertion holds for $k^{\prime}=n-1,n-2,\ldots k+1$. Then%
\begin{align*}
\delta_{n,k}  &  =\frac{q_{n-k}}{\prod\nolimits_{j=2}^{k}\sigma_{n-j}}%
\tilde{p}_{n-k-1}+\delta_{n,k+1}\\
&  \overset{\text{(\ref{defdeltan})}}{=}\frac{q_{n-k}}{\prod\nolimits_{j=2}%
^{k}\sigma_{n-j}}\tilde{p}_{n-k-1}+\frac{{\textstyle\prod\nolimits_{j=1}%
^{n-k-1}}\sigma_{j}} {{\textstyle\prod\nolimits_{j=1}^{n-2}} \sigma_{j}%
}\left(  {\displaystyle\prod\limits_{j=1}^{n-2-k}} \left(  1+q_{j+1}%
Q_{j}\right)  \right)  Q_{n-1-k}\\
&  \overset{\text{ind.}}{=}\prod\nolimits_{\ell=1}^{n-k-2}\left(  1+q_{\ell
+1}Q_{\ell}\right)  \left(  \frac{q_{n-k}}{\prod\nolimits_{j=2}^{k}
\sigma_{n-j}}+\frac{ {\textstyle\prod\nolimits_{j=1}^{n-k-1}} \sigma_{j}}{
{\textstyle\prod\nolimits_{j=1}^{n-2}} \sigma_{j}}Q_{n-1-k}\right) \\
&  =\frac{\prod\nolimits_{\ell=1}^{n-k-2}\left(  1+q_{\ell+1}Q_{\ell}\right)
}{\prod\nolimits_{j=2}^{k}\sigma_{n-j}}\left(  q_{n-k}+Q_{n-1-k}\right) \\
&  =\frac{\prod\nolimits_{\ell=1}^{n-k-2}\left(  1+q_{\ell+1}Q_{\ell}\right)
}{\prod\nolimits_{j=2}^{k}\sigma_{n-j}}Q_{n-k}\left(  1+q_{n-k}Q_{n-1-k}%
\right)  \sigma_{n-k}=\tilde{\delta}_{n,k}.
\end{align*}
Hence, (\ref{defdeltan}) is proved. We insert this into the right-hand side of
(\ref{ptildenrek}) and get%
\begin{align*}
\left(  \prod\nolimits_{\ell=1}^{n-2}\left(  1+q_{\ell+1}Q_{\ell}\right)
\right)  +\frac{q_{n}}{\sigma_{n-1}}\delta_{n,1}  &  =\left(  \prod
\nolimits_{\ell=1}^{n-2}\left(  1+q_{\ell+1}Q_{\ell}\right)  \right)  \left(
1+\sigma_{n-1}\frac{q_{n}}{\sigma_{n-1}}Q_{n-1}\right) \\
&  =\left(  \prod\nolimits_{\ell=1}^{n-2}\left(  1+q_{\ell+1}Q_{\ell}\right)
\right)  \left(  1+q_{n}Q_{n-1}\right) \\
&  =\prod\nolimits_{\ell=1}^{n-1}\left(  1+q_{\ell+1}Q_{\ell}\right)  =p_{n}.
\end{align*}

It remains to prove (\ref{doubleindb}). We insert (\ref{rekpn-1}) into the
right-hand side of (\ref{doubleindb}) and employ definition (\ref{defpntilde})%
\[
\det\mathbf{M}^{\left(  2n-1\right)  }=\left(  -1\right)  ^{n+1}\sum_{\ell
=1}^{n}\frac{q_{n+1-\ell}}{{\textstyle\prod\nolimits_{k=n+1-\ell}^{n-1}}
\sigma_{k}} \prod\nolimits_{k=1}^{n-\ell-1}\left(  1+q_{k+1}Q_{k}\right)  .
\]
We will prove%
\[
\det\mathbf{M}^{\left(  2n-1\right)  }=-\sigma_{n}Q_{n}\det\mathbf{M}^{\left(
2n\right)  }.
\]
For $1\leq k\leq n$, we introduce the partial sums%
\[
\lambda_{n,k}:=\sum_{\ell=k}^{n}\frac{q_{n+1-\ell}}{%
%TCIMACRO{\tprod \nolimits_{j=n+1-\ell}^{n-1}}%
%BeginExpansion
{\textstyle\prod\nolimits_{j=n+1-\ell}^{n-1}}
%EndExpansion
\sigma_{j}}\prod\nolimits_{j=1}^{n-\ell-1}\left(  1+q_{j+1}Q_{j}\right)
\]
so that $\det\mathbf{M}^{\left(  2n-1\right)  }=\left(  -1\right)
^{n+1}\lambda_{n,1}$. By induction for $k=n,n-1,\ldots,1$ we prove%
\begin{equation}
\lambda_{n,k}=\sigma_{n+1-k}\frac{\prod\nolimits_{j=1}^{n-k}\left(
1+q_{j+1}Q_{j}\right)  }{{\textstyle\prod\nolimits_{j=n+1-k}^{n-1}} \sigma
_{j}}Q_{n-k+1}. \label{deflambda}%
\end{equation}
We denote the right-hand side in (\ref{deflambda}) by $\tilde{\lambda}_{n,k}$
and prove $\lambda_{n,k}=\tilde{\lambda}_{n,k}$. For $k=n$, it holds
\[
\tilde{\lambda}_{n,n}=\frac{Q_{1}}{{\prod\nolimits_{\ell=2}^{n-1}}
\sigma_{\ell}}=\frac{q_{1}}{{\textstyle\prod\nolimits_{\ell=1}^{n-1}}
\sigma_{\ell}}=\lambda_{n,n}.
\]

We assume that (\ref{deflambda}) holds for $k^{\prime}=n,n-1,\ldots,k+1$.
Hence,%
\begin{align*}
\lambda_{n,k}  &  =\frac{q_{n+1-k}}{{\prod\nolimits_{j=n+1-k}^{n-1}}\sigma
_{j}}\prod\nolimits_{j=1}^{n-k-1}\left(  1+q_{j+1}Q_{j}\right)  +\lambda
_{n,k+1}\\
&  \overset{\text{(ind.)}}{=}\frac{q_{n+1-k}}{{\prod\nolimits_{j=n+1-k}^{n-1}%
}\sigma_{j}}\prod\nolimits_{j=1}^{n-k-1}\left(  1+q_{j+1}Q_{j}\right)
+\frac{{\prod\nolimits_{\ell=1}^{n-k-1}}\left(  1+q_{\ell+1}Q_{\ell}\right)
}{{\prod\nolimits_{\ell=n-k+1}^{n-1}}\sigma_{\ell}}Q_{n-k}\\
&  =\frac{\prod\nolimits_{j=1}^{n-k-1}\left(  1+q_{j+1}Q_{j}\right)  }%
{{\prod\nolimits_{j=n+1-k}^{n-1}}\sigma_{j}}\left(  q_{n+1-k}+Q_{n-k}\right)
\\
&  \overset{\text{(\ref{defamnqmntot})}}{=}\frac{\prod\nolimits_{j=1}%
^{n-k-1}\left(  1+q_{j+1}Q_{j-1}\right)  }{{\prod_{j=n+1-k}^{n-1}}\sigma_{j}%
}\left(  1+q_{n-k+1}Q_{n-k}\right)  Q_{n-k+1}\sigma_{n+1-k}\\
&  =\tilde{\lambda}_{n,k}%
\end{align*}
and this proves (\ref{deflambda}). We insert this into (\ref{deflambda}) to
obtain%
\[
\det\mathbf{M}^{\left(  2n-1\right)  }=\left(  -1\right)  ^{n+1}\lambda
_{n,1}=-\sigma_{n}Q_{n}\det\mathbf{M}^{\left(  2n\right)  }.
\]

\endproof

\appendix

\section{Some Basic Facts from Linear Algebra}

Note that the determinant of any $n\times n$ symmetric tridiagonal matrix
\[
\mathbf{W}_{n}:=%
\begin{pmatrix}
\gamma_{1} & \beta_{1} & 0 & 0\\
\beta_{1} & \gamma_{2} & \ddots & 0\\
0 & \ddots & \ddots & \beta_{n-1}\\
0 & 0 & \beta_{n-1} & \gamma_{n}%
\end{pmatrix}
\]
satisfies the three-term recursion
\begin{equation}
\det\mathbf{W}_{n}=\gamma_{n}\det\mathbf{W}_{n-1}-\beta_{n-1}^{2}%
\det\mathbf{W}_{n-2}. \label{detrek}%
\end{equation}
Let $\mathbf{W}_{n}^{\left(  i,j\right)  }$ denote the matrix which arises
when removing the $i$-th row and the $j$-th column. Then,%
\begin{equation}
\det\mathbf{W}_{n}^{\left(  i,n\right)  }=\left(
%TCIMACRO{\dprod \limits_{\ell=i}^{n-1}}%
%BeginExpansion
{\displaystyle\prod\limits_{\ell=i}^{n-1}}
%EndExpansion
\beta_{\ell}\right)  \det\mathbf{W}_{i-1}. \label{formulacof}%
\end{equation}

%----------------------------------------------------------------------------------------
%	Bibliography
%----------------------------------------------------------------------------------------

\bibliographystyle{abbrv}
\bibliography{bibfile}

\end{document}